\documentclass[10pt]{article}
\usepackage{amsmath,amsfonts,amssymb}
\usepackage{mathrsfs}
\usepackage{dsfont}
\usepackage[utf8x]{inputenc}
\usepackage[T1]{fontenc}
\usepackage{lmodern} \usepackage{ucs}
\usepackage{fullpage}

\usepackage{graphicx}
\usepackage{caption}
\usepackage{subcaption}

\usepackage{array}
\usepackage{multicol}
\usepackage{multirow}
\usepackage{color}
\usepackage{authblk}
\usepackage{esint}
\usepackage{enumitem,hyperref}

\usepackage{accents}

\usepackage{tikz}
\usetikzlibrary{decorations}
\usetikzlibrary{decorations.pathreplacing}
\usetikzlibrary{arrows}

\newtheorem{lemma}{Lemma}[section]
\newtheorem{theo}[lemma]{Theorem}
\newtheorem{rmk}[lemma]{Remark}
\newtheorem{proposition}[lemma]{Proposition}
\newtheorem{defin}[lemma]{Definition}

\newtheorem{coro}[lemma]{Corollary}

\makeatletter
\def\namedlabel#1#2{\begingroup
    #2%
    \def\@currentlabel{#2}%
    \phantomsection\label{#1}\endgroup
}
\makeatother

\makeatletter
\renewcommand*{\eqref}[1]{%
  \hyperref[{#1}]{\textup{\tagform@{\ref*{#1}}}}%
}
\makeatother

\newcommand{\QED}{\mbox{}\hfill \raisebox{-0.2pt}{\rule{5.6pt}{6pt}\rule{0pt}{0pt}} \medskip\par}

\newcommand{\R}{\mathbb{R}}

\def\E{\mathscr E}
\def\V{\mathscr V}
\def\pa{\partial}
\def\na{\nabla}

\newcommand{\ds}{\displaystyle}
\newcommand{\ud}{\, {\mathrm{d}}}

\newcommand{\eps}{\varepsilon}

\title{Constrained hydrodynamic flocking models in the limit of large attraction-repulsion interactions
}

\author[1]{Thierry~Goudon\thanks{ {\tt goudon@univ-cotedazur.fr}}}
\author[2]{Antoine Mellet\thanks{ {\tt mellet@umd.edu}}}
\date{}

\affil[1]{\small Universit\'e C\^ote d'Azur,   CNRS,  LJAD,  

Parc Valrose, F-06108 Nice, France}

\affil[2]{\small Department of Mathematics, University of Maryland, 

College Park MD 20742,  USA}

\begin{document}
\maketitle

\begin{abstract}
We study the collective dynamics of a population of particles/organisms subject to self-consistent attraction-repulsion interactions and an external velocity field. The starting point of our analysis is a  mean-field kinetic model and we investigate the singular limit corresponding to strong interaction forces. For well-prepared initial data, we show that the population asymptotically concentrates within a domain $\Omega(t)=\Omega_0+X(t)$ whose shape $\Omega_0$ is determined by the minimization of the interaction energy while the evolution of the domain’s center of mass $X(t)$ is determined by the external force field. 
In addition, we show that the internal flow of organisms within this moving domain is described by a  classical hydrodynamic model (the lake equation). The first part of our result relies only on the existence and uniqueness of minimizers for the interaction energy and holds for rather general interaction kernels. The second part is proved using a modulated energy method under more restrictive conditions on the nature of the interactions, and assuming that the limiting lake equation admits strong solutions.
\end{abstract}

\vspace*{.5cm}
{\small
\\[.4cm]

\noindent{\bf Math.~Subject Classification.} 
   35Q35, 
   35Q31, 
35Q83, 
35Q92 
}

\section{Introduction}

{\bf The mesoscopic model.}
The starting point of this paper is a mesoscopic description  of a large number of organisms (or particles/agents) moving under the combined actions of an attractive-repulsive force describing the interactions of the organisms with each others and a drag force due to an external velocity field.
The population of organisms is described by its  distribution 
 function in phase space $f_\eps(t,x,v)\geq 0$
 where $x$ denotes the position of the organisms and $v$ their velocity.
 The evolution of this function is described by the 
following Vlasov type kinetic equation:
 \begin{equation}
\label{VP1}
\partial_t f_\eps+v\cdot\nabla_x f_\eps - \eps^{-1}  \nabla_x \Phi_\eps \cdot\nabla_v f_\eps=\lambda\nabla_v\cdot 
((v-u_{ext})f_\eps).
\end{equation}
In the right hand side of \eqref{VP1}, $\lambda>0$ is a parameter (the strength of the drag force) and 
$(t,x)\in \mathbb R\times\mathbb R^N\mapsto u_{ext}(t,x)\in  \mathbb R^N$ is the given external velocity field.
Interactions between the organisms are described by the self-consistent potential $\Phi_\eps(t,x)$ defined by
\begin{equation}\label{VP2}
\Phi_\eps(t,x) = W  * \rho_\eps(t,x), \qquad   \rho_\eps(t,x)=\ds\int_{\R^N} f_\eps(t,x,v)\ud v ,
\end{equation}
where the kernel $W(x)$  combines the antagonistic effects of long-range attraction and short-range repulsion.
Equation \eqref{VP1} is set in the entire space $x\in \mathbb R^N$, $v\in\mathbb R^N$ and is supplemented with an initial data 
\begin{equation}
\label{alacharge}
f_\varepsilon(0,x,v)=f_\eps^{\mathrm{init}}(x,v)\geq 0,\qquad 
\ds\iint_{\R^N\times \R^N} f_\eps^{\mathrm{init}}(x,v)\ud v\ud x= m\in (0,\infty).
\end{equation}
The total  mass $m>0$ will be fixed throughout the paper.

\medskip

\noindent{\bf
Modeling issues and motivation.
}
Applications of system 
\eqref{VP1}-\eqref{VP2} include the description of the 
 behavior of a population of 
particles (in particular charged particles for applications in plasma physics) or organisms (cells, amoebae, fishes, birds etc), interacting with each
others via pairwise interactions that depend only on distance.
Such models have been used for 
describing collective behaviors in life sciences \cite{CuS2,CuS1,ToTu}, such as the 
flocking of birds \cite{Parrish},
the formation of ant trails \cite{Amo, Fourc}, the schooling of fish \cite{fish,Parr}, the swarms of bacteria \cite{bact,Xue}, 
  etc.
We refer to the kinetic model \eqref{VP1}-\eqref{VP2} as a mesoscopic model, as it can be derived as a mean-field limit of (microscopic) agent-based models \cite{Carr,Panf,Chayes,DegICM,Bert,HT}.

Interactions will be assumed to be repulsive at short distance and attractive at large distance.
It is well known that the interplay between these two antagonistic mechanisms leads to the formation of remarkable patterns
\cite{Barberis2016, Pain,Vic2} 
with potentially quite intricate  geometric structures.
The analysis of such phenomena has led to a very intensive research activity
\cite{Balag,Carr1,CDM,Choksi,Degond1,Lieb,Carr2, ShuTadmor,Simione}.
A recent overview with comments on relevant applications and impressive simulations can be found in \cite{Bailo}.
The study of these pattern structures largely relies on the analysis 
of the associated energy functional
\begin{equation}\label{energy}
\mathscr E[\rho]  : =   \iint_{\R^N\times\R^N} W(x-y) \rho (x) \rho(y) \ud x \ud y =  \int_{ \R^N} \Phi(x) \rho(x) \ud x
\end{equation}
 and on the properties of the energy minimizers 
within the set of 
non-negative measures with given mass $\int\ud \rho(x)=m$.
Existence, uniqueness and regularity properties of the global (and local) minimizers of the energy \eqref{energy}
are thoroughly discussed in 
\cite{Balag,Carr1,CDM,
CarrilloShu3d,CarrilloShu2d,CarrilloShu,
Choksi,
Lieb,Carr2, Ser_bk,
 ShuTadmor,Simione} and will play a key role in our analysis.

\medskip
\medskip

\noindent{\bf Hydrodynamic regime.}
The derivation of hydrodynamic models from \eqref{VP1}-\eqref{VP2} in the regime $\eps\to 0$ is inspired by the pioneering work of Brenier \cite{Brenier} for the Vlasov-Poisson equation on the torus ($x\in \mathbb T^N$, $-\Delta \Phi_\eps=\rho_\eps$):
imposing asymptotically a monokinetic particle distribution, 
one derives the incompressible Euler system (see also \cite{Mas} 
 for a refined asymptotic analysis). 
Motivated by applications in plasma physics, this question has been revisited 
in \cite{chiron}
when the equation is set in $\mathbb R^N$ with an external confining potential $\Phi_{\mathrm{ext}}$.
In that case, the hydrodynamic regime is determined by the opposite actions of this external potential and  the repulsive Coulomb forces.  
The effect of this competition is two-fold: 
it constrains the particles to a bounded domain $\Omega\subset \mathbb R^N$ and it prescribes the limiting macroscopic density (supported in that domain).
The flow of particles within this bounded domain is then described by the anelastic equation (or lake equation),  which is closely related to the Incompressible Euler system.
In the simplest case $\Phi_{\mathrm{ext}}(x)=\frac{|x|^2}{2}$,
the domain $\Omega$ is a ball, but it can be noticeably
more complicated in other cases. 
For example when $\Phi_{\mathrm{ext}}(x)=\Lambda x\cdot x$ with 
$\Lambda=\mathrm{diag}(1/\lambda_1^2,...,1/\lambda_N^2)\neq \mathbb I$, 
$\Omega$ has a surprising shape: it is an ellipsoid (as one would expect), but it does  not correspond to a level set of the external potential.
A related asymptotic analysis is carried out in \cite{serfaty,Ser_bk}, but starting   from a particle system, and combining the limit $\eps\to 0$ with the mean-field limit (number of particles goes to~$\infty$).

\medskip

The present paper  can be seen as a generalization of \cite{chiron} in two directions: First we work here with an attractive potential which is defined
self-consistently instead of being a given confining potential. Second, we consider  
 repulsive potentials that are not necessarily Coulombian.

\medskip

This paper was also motivated by some recent work of the second author and his collaborators on 
 another type of model that combines both attractive and repulsive forces: aggregation-diffusion equations. These are popular models for collective 
 dynamics that are very different in nature from \eqref{VP1}: They also involve a long range attractive force but they include a local repulsive force taken into account via a nonlinear diffusion term (overcrowding results in a pressure-like effect that drives the organisms away from each other).
In a series of work \cite{KMJW, KMW2,M24,MR25}, it was shown that singular limits of these equations lead to collective dynamics with sharp interfaces.
However, in that framework, the support of the limiting density changes with time and its evolution is determined by a free boundary problem of fluid type (such as Hele-Shaw flows with surface tension or mean-curvature flows depending on the regime of interest).
In the present work, we show that the competition between nonlocal attraction and nonlocal repulsion leads to a very different behavior and that the support of the density still moves, but only by translation of a fixed shape.

\medskip

Indeed, formally, the limit $f(t,x,v)  = \lim_{\eps\to 0} f_\eps(t,x,v)$ satisfies
$ \na_x \Phi \cdot \na_v f = 0 $. Multiplying this equality by $v$ and integrating with respect to $v$ yields
$$ \rho \na_x \Phi =0, \qquad \rho(t,x) = \int_{\R^N} f(t,x,v) \ud v , \quad \Phi = W*\rho.$$
This is the equation that characterizes the critical points of the  interaction energy \eqref{energy} under mass constraint. 
For well-prepared initial data, we will show that the limiting density must actually be a {\it global energy minimizer} of that energy. Under some assumptions on $W$, this prescribes the density profile  $\rho(t,x)$ for all time, up to a translation in space.

\medskip

\noindent{\bf Minimizers of the interaction energy.} 
With the mass $m>0$ fixed by the initial data \eqref{alacharge}, 
the following minimization problem thus plays a central role in our analysis (with $\E$ given by \eqref{energy}):
\begin{equation}\label{eq:min}
\mathscr E[\rho_0] = \min \left\{ \mathscr E[\rho]\, ; \, \rho\geq 0, \; \int_{\R^N} \rho \ud x =m\right\}.
\end{equation}
The first part of our analysis will be carried out for general kernels $W$ and  only requires \eqref{eq:min} to have a unique solution (up to translation). 
The second part of the paper requires additional properties of the energy $\E$ and its minimizers.
The full characterization of these minimizers is a delicate issue and for that reason, we restrict ourselves to kernels $W$ that are small perturbations of the classical repulsive power kernels with quadratic attraction:
\begin{equation}\label{kernW}
W(x) = - \alpha \frac{|x|^{p}}{p}  + \beta  \Lambda x\cdot x,\qquad 
\Lambda=\mathrm{diag}(1/\lambda^2_1,...,1/\lambda^2_N),   
\end{equation}
where $p\in(-N,-N+2]$ (with the usual convention that  $\frac{|x|^{p}}{p}$ becomes $ \ln|x|$ when $p=0$) and $\alpha,\beta>0$.
While the isotropic case $\Lambda=I$ is more classical, more general $\Lambda\neq I$
are relevant for many applications   involving  non-isotropic  interactions  (for fishes in a shallow river, for example, the vertical direction plays a different role from the horizontal ones).

When $W$ is given by \eqref{kernW} with $p\in(-N,2)$, it is well-known that there exists a unique measure $\rho_0(x)$,  solution of \eqref{eq:min} with center of mass at $0$ and supported on a compact set $\overline {\Omega_0}\subset \mathbb R^N$. 
With the further restriction $p\in(-N,-N+2]$, this minimizer is a bounded function   (and continuous when $p\in (-N,-N+2)$) - see for example \cite{CDM,FrMa}.
We note that $\rho_0$ is still absolutely continuous with respect to the Lebesgue measure when $p<-N+4$, while it is supported on a lower dimensional set when $p>-N+4$ - see Proposition \ref{facile} (it would be interesting to extend our analysis to this case, though the nature of the limiting equation changes significantly in that case).
We refer to \cite{Balag,FrMa} for further results on this issue.

\medskip

\noindent
{\bf The asymptotic model: informal statements.}
We now give an informal description of the main results of this paper, which characterize the asymptotic behavior of the solution $f_\eps(t,x,v)$ of \eqref{VP1} when $\eps\to0$.
We recall that $\rho_0(x)$ denotes the unique minimizer \eqref{eq:min} with center of mass at $0$, and we denote by $\overline {\Omega_0}$ its (compact) support.
We will prove the followings:
\begin{itemize}
\item[a)] Firstly, the density $\rho_\eps(t,x) = \int f_\eps(t,x,v)\, dv$ converges to $\rho(t,x)=\rho_0(x-X(t))$ where  
the center of mass $t\mapsto X(t)$ 
solves:
\begin{subequations}
     \begin{alignat}{1}\label{eq:XV1}
&\ds\frac{\ud}{\ud t}X (t) = V(t), \\[.3cm]
&\ds\frac{\ud}{\ud t} V(t) = \ds\frac{\lambda}{m} \int_{\R^N}  \rho_0(x-X(t)) ( u_{ext} (t,x)- V(t) )\ud x.\label{eq:XV2}
\end{alignat}
\end{subequations}
When endowed with  initial conditions
\begin{equation}\label{eq:initas}
X(0)= X^{\mathrm{init}}, \qquad V(0)=V^{\mathrm{init}}, \end{equation}
 this system admits a  unique solution $(X,V)$, see Lemma~\ref{Lem:edo}.

The limit density $\rho(t,\cdot)$ is thus  supported on the moving set $\Omega(t) = \Omega_0 + X(t)$.
This result establishes the swarming behavior (or phase separation) in the limit $\eps\to0$: it shows  the formation of a moving interface $\pa\Omega(t)$ separating the vacuum region from the region of positive density. 
A more precise statement is given in Theorem \ref{th1} and will be proved using basic a priori energy estimates and assumptions on the initial data which ensure that the limit density $\rho(t,x)$ minimizes the interaction energy $\mathscr E$ for all $t>0$. 
This result is quite general and holds for very general kernel $x\mapsto W(x)$.

\item[b)] Our second result shows that the macroscopic flux $j_\eps(t,x)=\int_{\mathbb R^N} vf_\eps(t,x,v)\ud v$ converges, when $\eps\to0$, to $j(t,x)=\rho(t,x) \V(t,x) = \rho_0(x-X(t)) \V(t,x)$ where the velocity field $\V(t,x)$ is solution of 
\begin{equation}\label{eq:asymptotic}
\begin{cases}
\pa_t \rho + \na_x \cdot (\rho \V)=0 & \mbox{ in } \Omega(t), \\
\pa_t (\rho \V) + \na_x(\rho \V\otimes \V) + \rho \na_x P = \lambda \rho(u_{ext}-\V) & \mbox{ in } \Omega(t), \\
\rho \V\cdot \nu  = \rho V\cdot \nu   & \mbox{ on } \pa \Omega(t),
\end{cases}
\end{equation}
with $\nu$ standing for the outward unit vector on $\partial\Omega(t)$ (and $V(t)$ is the velocity appearing in \eqref{eq:XV1}-\eqref{eq:XV2}).
The equation is supplemented with an initial condition 
 \begin{equation}\label{eq:initas2}
 \V(0,x)=\V^{\mathrm{init}}(x)
\end{equation}
where the functions 
$\rho^{\mathrm{init}}(x):=\rho_0(x- X^{\mathrm{init}})$ and $ \V^{\mathrm{init}}(x)$ satisfy the compatibility condition
\begin{equation}\label{eq:initcomp}
\begin{cases}
\ds V^{\mathrm{init}} = \frac{1}{m} \int_{\R^N} \rho^{\mathrm{init}}(x) \V^{\mathrm{init}}(x)\ud x, \\[5pt]
\ds \na_x \cdot (\rho^{\mathrm{init}}[ \V^{\mathrm{init}}-V^{\mathrm{init}}]) = 0 \mbox{ in } \Omega^{\mathrm{init}}, \qquad \rho^{\mathrm{init}}[ \V^{\mathrm{init}}-V^{\mathrm{init}}] \cdot \nu  = 0 \mbox{ on } \pa\Omega^{\mathrm{init}}.
\end{cases}
\end{equation}
System \eqref{eq:asymptotic} is called the anelastic system, or ``lake equations'', and it is set here  on the  moving domain $\Omega(t)$ (when $\Omega(t)=\Omega_0$ is fixed and $x\mapsto \rho_0(x)$ is constant in $\Omega_0$, then \eqref{eq:asymptotic} is the standard incompressible Euler system).
In the context of collective motion, we note that our first result above identifies the limiting density, while this second result tracks the motion (or flow) of the organisms/particles inside $\Omega(t)$.

 We refer the reader to 
 \cite{Lake1,Lake2}
 for an introduction to the lake equation and for some existence results (for strong solutions), when the density $\rho$ is assumed to be bounded from below away from zero.
 We also refer to \cite{Mas2} for a rigorous derivation of this system from the compressible Navier-Stokes system,
and to \cite{DuSe} for a derivation  as a mean-field limit for Ginzburg-Landau vortices.

The rigorous justification of this second result (a precise statement of which can be found in Theorem \ref{th2}) will require additional assumptions on the minimization problem \eqref{eq:min}.
The proof relies on a modulated energy argument, in the spirit of \cite{Brenier}, which requires in particular the existence of a strong solution of the limiting equation \eqref{eq:asymptotic} - something that is typically only achievable on a finite time interval. Our approach is similar to the method that was used in \cite{chiron} (when $p=-N+2$) and \cite{serfaty} (when $p\in(-N,-N+2]$) to study a related problem, in which the pairwise attractive interactions are replaced by a given confining potential (so that the particles are constrained to a fixed domain). 
Compared to these papers, our definition of the modulated energy will need to account for the motion of the center of mass.
\end{itemize}

We can summarize our results as follows:
The energy minimization problem \eqref{eq:min} identifies  the profile of the limiting density $\rho(t,x)$ and its support $\Omega(t)$ up to translation.
The time evolution of this density is reduced to the motion of its center of mass, which is driven by the external velocity field $u_{ext}(t,x)$ via \eqref{eq:XV1}-\eqref{eq:XV2}.
Inside the domain, under appropriate regularity assumptions, the flow of the organisms/particles is described by the lake equation \eqref{eq:asymptotic}.

\medskip

{\bf Outline of the paper.}
Section~\ref{Prel} will make precise the functional framework: we collect there the  main assumptions on the parameters of the model 
($p$, $N$, $u_{ext}$...). We will then state precisely our main results (Theorems \ref{th1} and  \ref{th2}).
In Section~\ref{Sec:Formal} we present the formal derivation of the asymptotic system, and point out the difficulties in making the derivation rigorous.
Section~\ref{sec:min} includes several important statements regarding  the minimization problem \eqref{eq:min} as well as key coercivity estimates for the interaction energy $\mathscr E$.  In Section~\ref{sec:continuity}, we establish the convergence of the density (proof of Theorem \ref{th1}) and in Section~\ref{Sec:flux}, we turn to the modulated energy analysis and the convergence of the flux (proof of Theorem \ref{th2}).
Finally, Section~\ref{sec:ell} provides some explicit computations concerning the solution $\rho_0(x)$ of the minimization problem for non-isotropic kernel $\Lambda\neq I$ allowing us to rigorously justify some of our key assumptions in some particular cases.

\section{Preliminaries and Main Results}
\label{Prel}

\subsection{Notations and Main Assumptions}
We recall that the external velocity field $u_{ext}:\mathbb R\times\mathbb R^N\to \mathbb R^N$ is given, and we assume that it satisfies
\begin{equation}\label{hyp_uext}
u_{ext}\in  C^0(\mathbb R\times \mathbb R^N)\cap L^\infty(\mathbb R\times \mathbb R^N)\cap L^\infty(\mathbb R; W^{1,\infty}(\mathbb R^N)).
\end{equation}
 
\medskip

Next, we denote by $\mathcal M^1(\R^N)$ the set of Radon measures in $\R^N$ and $\mathcal M^1_+(\R^N)$ the set of non-negative Radon measures in $\R^N$ and we consider the following minimization problem (with $\E$ defined by \eqref{energy}):
\begin{equation}\label{eq:min0}
\mathscr E_m = \inf \left\{ \mathscr E[\rho] \, ;\, \rho \in \mathcal M^1_+, \; \int_{\mathbb R^n} \ud\rho(x)=m\right\}.
\end{equation}
Since the energy $\mathscr E$ is  invariant by translation ($\mathscr E[\rho(\cdot-X)] = \mathscr E[\rho]$ for all $X\in \R^N$), it is natural 
to use the position of the center of mass $\frac 1 m\int_{\mathbb R^N} x\ud \rho(x)$ as a selection criterion:
a key property needed for our analysis is the uniqueness of the minimizer $\rho_0$ with center of mass at $x=0$. In fact, our first result will be proved under the following general assumptions:
\begin{itemize}
\item[{(\bf H1a)}] The interaction kernel $W:\R^N \to (-\infty,+\infty]$ is locally integrable, lower semi-continuous, bounded below and symmetric ($W(-x)=W(x)$).
\item[{(\bf H1b)}] There exists a unique $\rho_0 \in \mathcal M^1_+$ minimizer of \eqref{eq:min0} with center of mass at $0$.
\end{itemize}

Assumption {(\bf H1a)} implies in particular that $\rho\mapsto \mathscr E[\rho]$ is lower-semicontinuous with respect to weak convergence of measures (see for instance \cite[Lemma 2.2]{Simione}).
The existence of a global minimizers then follows under mild assumptions on the behavior of $W$ at $\infty$ \cite{Simione} (see also \cite{Choksi} for results with power law kernels). 
The uniqueness of this minimizer, however, is not guaranteed  and has only been proved under further structural assumptions on $W$ (see \cite{CarrilloShu,DYY,Lopes} and references therein).
For our second result, we restrict ourselves to simple $W$ for which {(\bf H1b)} is known to hold:  
\begin{itemize}
\item[{(\bf H2a)}] $W$ is given by 
\begin{equation*}\label{eq:WW}  
W(x) = - \alpha \frac{|x|^{p}}{p}  + \beta  \Lambda x\cdot x+w(x),\qquad 
\Lambda=\mathrm{diag}(1/\lambda^2_1,...,1/\lambda^2_N), 
\end{equation*}
with $\alpha,\beta>0$, $p\in (-N,2)$ and  where $w(x)$ satisfies:
\begin{equation*}\label{eq:H1}
\begin{array}{c}
\mbox{ $w:\R^N\to (-\infty,+\infty]$ is locally integrable, lower semi-continuous, } \\
 \mbox{bounded below  and symmetric ($w(x) = w(-x)$)}  
\end{array}
\end{equation*}
and  is such that
\begin{equation*}\label{compinf}
\ds\lim_{|x|\to \infty}  \frac{w(x) }{|x|^2} =0.
\end{equation*}
\item[{(\bf H2b)}] In addition, we assume that 
\begin{equation*}\label{eq:conditionsurp}
p\in(-N,-N+2]\quad  \mbox{ and } \quad 
\mathscr F[w](\xi) \geq -\kappa \mathscr F\left[- \alpha \frac{|x|^{p}}{p} \right](\xi) \; \forall \xi\in \mathbb R^N \mbox{ for some  $\kappa\in(0,1)$}
\end{equation*}
where $\mathscr F$ stands for the Fourier transform.
\item[{(\bf H2c)}]  Finally, $w$ satisfies
\begin{equation*}\label{eq:wreg}
\xi\mapsto |\xi|^{(p+N)/2}\mathscr F(\partial^k w)(\xi)
\in L^2(\mathbb R^N) 
\qquad\mbox{ for all $k\in \mathbb N^N$ such that $1\leq  |k|\leq 2$}.
\end{equation*}
\end{itemize}
 
We  note that {(\bf H2a)} immediately implies {(\bf H1a)}, while we will see that {(\bf H2b)} implies {(\bf H1b)}. So these assumptions {(\bf H2)}
identify a particular class of kernels for which {(\bf H1)} holds.
The last condition {(\bf H2c)} is a technical assumption that will be used in the proof of Theorem \ref{th2}. Note that  
{(\bf H2b)}
can be viewed as a smallness condition on $w$ but that there is no smallness requirement in the regularity assumption {(\bf H2c)}. 

We recall that when $\Lambda = \mathbb I$ and $w=0$, the properties of the unique minimizer $\rho_0$ are well known: we have
\begin{proposition}[Case $\Lambda = \mathbb I$, $w=0$]\label{facile}
Assume that $W$ is given by 
{(\bf H2a)}
with $p\in (-N,2)$, $\Lambda = \mathbb I$ and $w=0$. Then for all $m>0$, there exists a 
unique 
minimizer $\rho_0\in \mathcal M^1_+$  for \eqref{eq:min0} with center of mass at $0$. Moreover, the following assertions hold:
\begin{itemize}
    \item when $-N<p < -N+2$, $\rho_0$ is a continuous function in $\R^N$ and is supported on $\overline {B_{R_0}}=\{x\in \R^N\,;\, |x|\leq R_0\}$,
    \item when $p=-N+2$, $\rho_0$ is a function in $L^\infty(\R^N)$ and is given by 
    $\rho_0(x)= c_0\chi_{B_{R_0}}(x)$;
    \item when $-N+2<p<\min\{2,-N+4 \}$, $\rho_0$ is a function in $L^1(\R^N)$ supported on a ball $\overline{B_{R_0}}$.
    \item when $-N+4\leq p<2$, the minimizer $\rho_0$ is a measure supported on a sphere $\pa B_{R_0}$.
\end{itemize}
In all cases, the radius $R_0$ depends on $\alpha$, $p$, $\beta$ and $m$.
\end{proposition}
We refer to \cite{FrMa} for a synthesis of the literature on this problem (in particular the work in \cite{Imbert,CaVa,CarrilloShu}) and for various explicit formulas.
Extension of these results can be proved when $\Lambda\neq \mathbb I$, although the radial symmetry is broken in that case. 
When we include the perturbation $w$, condition 
{(\bf H2a)}
is enough to guarantee  the existence of a compactly supported global minimizer $\rho_0$ while 
condition 
{(\bf H2b)}
yields  the uniqueness. This will be further discussed in Section \ref{sec:min1}.

\medskip

The potential  $\Phi_0=W*\rho_0$ associated to the minimizer $\rho_0$ also plays an important role in the proofs. 
An important characterization of the minimizers (see \cite{Balag}) is the fact that 
there exists a constant $A_0$ such that
$$ \Phi_0(x)\geq A_0  \mbox{ in } \mathbb R^N \mbox{ and  } \Phi_0(x)= A_0 \mbox{ on }\mathrm{supp}(\rho_0)$$
(with $mA_0=\mathscr E_m$).
In the framework of Proposition \ref{facile}, we can then show the following property:
\begin{equation}\label{eq:conditionsurphi}
\begin{split}
& \mbox{For any Lipschitz vector field $\mathscr V_0:\mathbb R^N\to \mathbb R^N$ such that $\mathscr V_0\cdot \nu =0 $ on $\pa\Omega_0$,}\\
& \mbox{there exists $C_0>0$ such that }
|\mathscr V_0(x)\cdot \nabla\Phi_0(x)|\leq C_0(\Phi_0(x)-A_0) \qquad \forall x\in \mathbb R^N.
\end{split}
\end{equation}
This property turns out to be crucial for the asymptotic analysis and to derive the lake equation.
It holds under our assumptions {(\bf H2a)}-{(\bf H2b)} at least when $w=0$.
This was proved in \cite{chiron} (for the case $p=-N+2$) and 
in \cite{serfaty} (for the case $p\in (-N,N+2]$). 
In Section~\ref{sec:ext} we present some computations that establish this inequality when $w=0$.
We will not attempt to generalize these computations to the case $w\neq 0$,
and we will instead add the following assumption:
\begin{itemize}
\item[{(\bf H2d)}] We further assume that the global minimizer $\rho_0$ of $\mathcal E$ is such that \eqref{eq:conditionsurphi} holds with $\Phi_0=W*\rho_0$.
\end{itemize}

\subsection{Energy Functional for the Vlasov Equation  \eqref{VP1}}
\label{Sec:enfunct}
Weak solutions of \eqref{VP1} can be  defined  as functions $f \in L^\infty(0,T;L^1(\R^N\times\R^N))$ that satisfy the following weak formulation 
$$\begin{array}{l}
\ds\int_0^T \iint_{\R^N\times\R^N} f [\pa_t\psi+ v\cdot\nabla_x\psi - \eps^{-1}  \nabla_x \Phi \cdot\nabla_v\psi -\lambda 
(v-u_{ext})\cdot  \nabla_v \psi ] \ud x \ud v \ud t
\\
\hspace*{5cm}=- \ds\iint_{\R^N\times\R^N} f^{\mathrm{init}} (x,v)\psi(0,x,v) \ud x\ud v  
\end{array}$$
for any smooth  test function $\psi$,  compactly supported in $[0,\infty)\times\mathbb R^N\times\mathbb R^N$.
These solutions should also satisfy the mass conservation property
\begin{equation}\label{mass_cons}
\ds\iint_{\R^N\times\R^N}  f_\eps(t,x,v)\ud x\ud v=\ds\iint_{\R^N\times\R^N}  f^{\mathrm{init}}_\eps(x,v)\ud v\ud x
=m,
\end{equation}
as well as the natural energy inequality associated to \eqref{VP1}:
Given $f:\mathbb R^N\times\mathbb R^N\to [0,\infty)$, the total energy is defined by 
\begin{equation}\label{def_H}
\mathscr H_\eps[f] = \iint_{\R^N\times\R^N} \frac{|v|^2}{2} f(x,v)\ud x\ud v + \frac{1}{2\eps} \left(
\mathscr E[\rho] - \mathscr E_m\right)\geq 0, \qquad \rho(x) =\int_{\R^N} f(x,v)\ud v.
\end{equation}
We recognize the sum of the usual  kinetic energy and the (normalized) potential energy of the system.
Solutions of \eqref{VP1} then satisfy
\begin{align}
\ds\frac{\ud}{\ud t} 
\mathscr H_\eps[f_\eps(t,\cdot)] 
& = - \lambda\iint_{\R^N\times\R^N}  v\cdot(v-u_{ext})f_\eps\ud v\ud  x \label{evol_H} \\ 
& \leq -\frac{\lambda}{2} \iint_{\R^N\times\R^N}  |v|^2 f_\eps\ud v\ud  x + \frac{\lambda}{2} \int_{\R^N}  |u_{ext}|^2 \rho_\eps\ud  x.\nonumber 
\end{align}
When $u_{ext}$ satisfies \eqref{hyp_uext} this inequality implies that the energy $\mathscr H_\eps[f_\eps(t)]  $ is bounded uniformly with respect to $\eps>0$ and $t\in [0,T]$ if it is initially bounded (see Proposition~\ref{prop:bounds}).
It follows that $\mathscr E[\rho_\eps(t)] \to \mathscr E_m$
This is the key to identifying the limiting density $\lim_{\eps\to 0} \rho_\eps$
(using ({\bf H1a}) and ({\bf H1b})).

We sketch the construction of weak solutions of  \eqref{VP1} satisfying \eqref{mass_cons} and \eqref{evol_H}  in Section~\ref{WeakSol}.


\medskip

Finally, given solution $f_\eps(t,x,v)$ solution of \eqref{VP1}-\eqref{VP2}, we define
the macroscopic density and flux
\begin{equation}\label{eq:rhou}
\rho_\eps (t,x)= \int_{\R^N} f_\eps(t,x,v)\ud v, \qquad j_\eps(t,x) = \int_{\R^N} v f_\eps(t,x,v)\ud v,
\end{equation}
as well as the center of mass and average velocity
\begin{equation}\label{eq:XVepsdef}
X_\eps(t) = \frac{1}{m}\ds \iint_{\R^N\times\R^N } x f_\eps(t,x,v)\ud v\ud x,\qquad V_\eps (t) = 
\frac{1}{m} \ds\iint_{\R^N\times\R^N}  vf_\eps (t,x,v)\ud v \ud x.
\end{equation}

\subsection{Limit $\eps\to 0$: Main Theorems}
Our first result will be proved for initial data $f_\eps^{\mathrm{init}}:\mathbb R^N\times\mathbb R^N\to [0,\infty)$
that satisfy a uniform bound on the energy:
 \begin{equation}
 \label{bornener}
 \ds\iint_{\R^N\times\R^N}f_\eps^{\mathrm{init}}\ud v\ud x=m,\qquad 
  \sup_{0 < \eps < 1} \left\{
  \ds\iint_{\R^N\times\R^N} \left( \frac{ |v| ^2}2 + \frac{|x|^2}{2}\right)\ds\ f_\eps^{\mathrm{init}} \ud v\ud x 
 + \ds\frac{\mathscr E[\rho_\eps^{\mathrm{init}}]- \mathscr E_m}{2\eps}\right\}
 <\infty.
\end{equation}
This condition implies in particular that the initial density $\rho_\eps^{\mathrm{init}}$ is a small perturbation of an energy minimizer $\rho_0(x-X^{\mathrm{init}})$ (with  $X^{\mathrm{init}}\in \R^N$). 

The second result requires an initial data that is  a small perturbation of the macroscopic equilibrium obtained as a minimizer of the energy, that is 
\begin{equation}\label{eq:rhoint}
 f_\eps^{\mathrm{init}} (x,v)\sim  \rho_0(x-X^{\mathrm{init}}) \delta(v-\V^{\mathrm{init}}(x))
\end{equation}
where  $\V^{\mathrm{init}} $ is a 
smooth vector field in $\Omega^{\mathrm{init}}=\Omega_0+X^{\mathrm{init}}$
satisfying the compatibility conditions \eqref{eq:initcomp}.
To be more specific, we will require  $f_\eps^{\mathrm{init}}$ to satisfy 
\begin{align}\label{eq:initass}
& \lim_{\eps\to 0}  \left[ \ds\frac{1}{2} \iint_{\R^N \times\R^N}  |v-\mathscr V^{\mathrm{init}}(x)|^2 f_\eps^{\mathrm{init}}(x,v)\ud v\ud x \right.\nonumber \\
& \qquad\qquad\qquad  \left.  + 
 \frac{1}{2\eps} \Big( \mathscr E[\rho_\eps^{\mathrm{init}}]-\mathscr E_m + |X_{\eps}^{\mathrm{init}}-X^{\mathrm{init}}|^2 
 +|V_\eps^{\mathrm{init}}-V^{\mathrm{init}})|^2\Big)\right] =0 
\end{align}
 where the initial center of mass and average velocity are defined (according to \eqref{eq:XVepsdef}) by:
$$
X_\eps^{\mathrm{init}}=\ds\frac1m\ds\iint_{\R^N\times\R^N} xf_\eps^{\mathrm{init}}(x,v)\ud v \ud x, \qquad V_\eps^{\mathrm{init}}=\ds\frac1m\ds\iint_{\R^N\times \R^N} v f_\eps^{\mathrm{init}}(x,v)\ud v \ud x.$$
We can now state our main results.
First we have
\begin{theo}
\label{th1} 
Let $u_{ext}$ satisfy \eqref{hyp_uext} 
and let $W$ be such that ({\bf H1a}) and ({\bf H1b}) hold.
Consider initial conditions satisfying the uniform energy bound \eqref{bornener} and such that 
$$\lim_{\eps\to 0} (X_\eps^{\mathrm{init}},V_\eps^{\mathrm{init}})  = (X^{\mathrm{init}},V^{\mathrm{init}}).
$$
Let $f_\eps(t,x,v)$ be the associated weak solution  of the Vlasov equation 
\eqref{VP1}. Then, for any $0<T<\infty$, we have:
\begin{itemize}
\item[(i)] the center of mass and velocity $t\mapsto (X_\eps(t), V_\eps(t))$ defined by \eqref{eq:XVepsdef} converge uniformly in $[0,T]$ to $t\mapsto (X(t),V(t))$ solution of  \eqref{eq:XV1}-\eqref{eq:XV2} with initial conditions  $(X^{\mathrm{init}},V^{\mathrm{init}})$.
\item[(ii)] the density $\rho_\eps(t,x)$ converges to $\rho(t,x)=\rho_0(x-X(t))$ in 
$C^0([0,T];\mathscr M^1(\mathbb R^N)-\text{weak}-\star)$.

\item[(iii)] Up a subsequence, $j_\eps(t,x)$ converges to $j(t,x)$ in $\mathscr M^1([0,T]\times
  \mathbb R^N)$ weakly-$\star$ which satisfies
    \begin{equation}\label{lim_mc}
\partial_t \rho+\nabla_x\cdot j=0 \qquad\mbox{ in } [0,T]\times\mathbb R^N.
\end{equation}
Moreover,  $j$ is absolutely continuous with respect to $\rho$. 
\end{itemize}
\end{theo}

We point out that Theorem~\ref{th1} holds in a very general framework. In particular it does not assume any regularity on the minimizer $\rho_0$ as long at it exists and is unique (up to translation). 
This holds when $W$ is given by 
{(\bf H2a)} with a repulsive power in the whole range  $p\in (-N,2)$ which includes values of $p$ for which the minimizer is a Radon measure rather than a function.

The proof of this result is given in Section \ref{sec:continuity} and follows from the following steps:
\begin{enumerate}
\item First we prove that $\rho_\eps(t,x)$ and $j_\eps(t,x)$ have limits $\widetilde \rho(t,x)$ and $j(t,x)$ (up to subsequences) which satisfy the continuity equation $\pa_t \widetilde \rho + \na_x\cdot j=0$
and have the form $\widetilde \rho(t,x) = \rho_0(t-\widetilde X(t))$ for some $\widetilde X(t)\in \R^N$.
\item We will then show that $\widetilde X(t)=X(t)$, the solution of \eqref{eq:XVlim}, and deduce that $\widetilde \rho(t,x) = \rho(t,x) =\rho_0(t- X(t))$. 
\end{enumerate}
\medskip

Next, with strengthened assumptions, we can prove the convergence of the flux toward a solution of  the lake equation (provided such a solution exists):
\begin{theo}
\label{th2} 
Let $W(x)$ be given by ({\bf H2a})-({\bf H2d})
and assume that the initial conditions satisfies \eqref{bornener} and \eqref{eq:initass}.
Assume further that:
\begin{itemize}
\item[] The system  \eqref{eq:XV1}-\eqref{eq:XV2}-\eqref{eq:asymptotic} with initial conditions \eqref{eq:initas} has a strong solution $(\rho(t,x),X(t),V(t),\V(t,x))$ defined on a time interval $ [ 0, T ] $ in the sense of Definition \ref{def:strong}. 
and the function $\V(t,x)$ has a Lipschitz extension (still denoted by $\V$) to $[0,T]\times \R^N$.
\end{itemize}
Then, the whole sequence $j_\eps(t,x)$ converges to $\rho(t,x) \mathscr V(t,x)$ in $\mathscr M^1([0,T]\times
  \mathbb R^N)$.    Furthermore, we have 
$$ \lim_{\eps\to 0}  \ds\sup_{0\leq t\leq T}  \iint_{\R^N \times\R^N} \frac{|v-\mathscr V(t,x)|^2}{2}  f_\eps(t,x,v)\ud v\ud x =0.
$$ 
 \end{theo}

\noindent
The proof of this second theorem is more delicate and relies on the modulated energy method introduced in \cite{Brenier}.
We note that the uniqueness of the limits imply that the whole sequences $\rho_\eps(t,x)$ and $j_\eps(t,x)$ converge. Furthermore, $f_\eps(t,x,v)$ converges to a function $f(t,x,v)$ satisfying
\[
\mathrm{supp}(f)\subset [0,T]\times
\overline {\Omega(t)}\times\R^N.
\]
Extending Theorem \ref{th2} to more general interaction kernels $W$ requires some additional work. In particular: 
\begin{itemize}
\item The definition of the modulated energy \eqref{eq:H} and the resulting inequality  (Proposition \ref{Pr:estGron}) would require a different formulation for non-quadratic attractive kernels.
\item The proof uses in a crucial way the estimate \eqref{eq:conditionsurphi} which we will only prove when $w=0$. 
Nonetheless, this estimate is very natural in view of the expected regularity and non-degeneracy of the solution of the associated obstacle problem that characterize the minimizer $\rho_0$ (see \cite{CDM}). It is expected to hold under some regularity and smallness assumptions on $w$, and can probably be extended to other powers as long as the minimizer $\rho_0$ is an $L^1$ function, that is for  $p\in (-N,-N+4)$.
\item When $W$ is given by {(\bf H2a)} 
with $p\in [-N+4,2)$, the minimizer $\rho_0$ is a singular measure supported on a lower-dimensional set.
It would be interesting to understand the asymptotic dynamics and the meaning of the 
 macroscopic equation in that case.
 \item Finally, we note that the present analysis can likely be extended,
 using the approach of \cite{serfaty} and the generalization of the modulated energy \cite{Ser_bk}, 
  to   study the limit of the corresponding particle system.   
\end{itemize}

\section{The Asymptotic System}
\label{Sec:Formal}

\subsection{Formal Derivation}
In this section, we formally derive the asymptotic model in order  to understand 
how the system \eqref{eq:XV1}-\eqref{eq:XV2}-\eqref{eq:initas}-\eqref{eq:asymptotic} emerges in the limit $\eps\to 0$. 
Some of the computations presented in this section will also be useful in the proofs of the main theorems.

First, we integrate \eqref{VP1} with respect to $v$ to get the continuity equation
\begin{equation}\label{eq:rhoeps}
\partial_t \rho_\eps + \nabla_x\cdot j_\eps=0.
\end{equation}
Multiplying \eqref{VP1} by $v$ and then integrating with respect to $v$ gives the momentum equation
\begin{equation}\label{eq:jeps}
\partial_t j_\eps + \mathrm{div}_x \mathbb P_\eps+  \rho_\eps \nabla_x\frac{ \Phi_\eps}{\eps} =\lambda( \rho_\eps u_{ext}-j_\eps),
\end{equation}
where we have set 
\[\mathbb P_\eps(t,x)=\int_{\R^N} v\otimes v f_\eps(t,x,v)\ud v.\]
We now explain how to pass to the limit (formally) in \eqref{eq:rhoeps}-\eqref{eq:jeps}.
\\

\noindent{\bf The density.}
The bound on the energy $\mathscr H_\eps[f_\eps(t)]$, see \eqref{def_H} and \eqref{evol_H},together with the lower semi-continuity of $\mathscr E$ (Assumption ({\bf H1a})) implies that the limiting density $\rho(t,x)$ satisfies
$\mathscr E[\rho(t,\cdot)] = \mathscr E_m $ and so $x\mapsto \rho(t,x)$ is a (global) minimizer of the interaction energy $\mathscr E$, with mass constraint,  for all $t>0$.
The uniqueness of these minimizers (Assumption ({\bf H1b})) thus implies that
\begin{equation}\label{eq:rhol}
\rho(t,x)  = \rho_0(x-X(t))
\end{equation}
for some $X(t)\in \R^N$. 
We note that $X(t)$ is the center of mass, defined by
$$ X(t) = \frac{1}{m}\ds \int_{\R^N } x \rho(t,x) \ud x$$
and thus $\mathrm{supp}(\rho(t,\cdot))=\Omega(t)=\Omega_0+X(t)$ is translated from $\mathrm{supp}(\rho_0)=\Omega_0$.
\\

\noindent{\bf The continuity equation.}
Next, we assume that $j_\eps(t,x) \to j(t,x)$. 
Passing to the limit in the continuity equation \eqref{eq:rhoeps} then yields
\begin{equation}\label{eq:rho1}
\pa_t \rho + \na_x \cdot j=0 \qquad \mbox{ in } (0,\infty)\times\R^N.
\end{equation}
We define the macroscopic velocity
$$
V  (t) = 
\frac{1}{m} \ds\int_{\R^N}  j (t,x) \ud x
$$
and we note that multiplying \eqref{eq:rho1} by $x$ and integrating implies
$$ X'(t) = V(t)$$
which is \eqref{eq:XV1}.

Differentiating  the relation \eqref{eq:rhol}, we find $\pa_t \rho = -\na_x \rho \cdot X'(t) = -\na_x \rho \cdot V(t) $ so \eqref{eq:rho1} can also be written as
\begin{equation}\label{eq:rho2}
\na_x \cdot (j-\rho V)=0 \qquad \mbox{ in } (0,\infty)\times\R^N.
\end{equation}
The fact that this continuity equation is satisfied in $\R^N$, and not just on the support of $\rho$, is important since it means that it encodes the null-flux boundary conditions on the boundary of that support:
since $\rho(t,\cdot)$ is supported in $\Omega(t)=\Omega_0+X(t)$, this equality is in fact equivalent to 
\begin{equation}\label{bc_fluid}
\na_x \cdot (j- \rho V(t)) = 0 \qquad \mbox{ in }\Omega (t), \qquad (j-\rho V(t))\cdot \nu = 0 \qquad \mbox{ on }\pa \Omega (t).
\end{equation}

 \medskip

\noindent{\bf The velocity equation.}
We define $\V(t,x)$ such that $j(t,x) = \rho(t,x) \V(t,x)$ (this implicitly assumes that $j$ is absolutely continuous with respect to $\rho$).  
Using \eqref{eq:rho1}, \eqref{eq:rho2}, we get
\begin{equation}\label{eq:rho1bis}
\pa_t \rho + \na_x \cdot (\rho \V)
=\na_x \cdot (\rho (\V-V))=
0 \qquad \mbox{ in } (0,\infty)\times\R^N.
\end{equation}
We turn  to the passage  to the limit in \eqref{eq:jeps}.
Since $j_\eps$ converges to $\rho \mathscr \V$, we expect $\mathbb P_\eps$ to converge to $\rho \mathscr V\otimes V$.
This is in line  with the fact that we expect the initial monokinetic profile \eqref{eq:rhoint} to propagate to the solution, that is
$f_\eps(t,x,v)\sim \rho(t,x)\delta(v-\V(t,x))$
(The modulated energy is designed to keep track of this behavior).

Furthermore, we will see that the term  $ \rho_\eps \nabla_x\frac{ \Phi_\eps}{\eps}$ is bounded and admits a limit. This limit is of the form $\rho\nabla_xP$, where 
$P$ can be interpreted as a Lagrange multiplier associated to the constraint \eqref{eq:rho1bis}.
With these considerations, \eqref{eq:jeps} leads to
$$
\partial_t (\rho \mathscr V)+\na_x (\rho \V\otimes \mathscr V) +\rho \nabla_xP=\lambda\rho (u_{ext}-\mathscr V) \qquad \mbox{ in } (0,\infty)\times\Omega(t).
$$
This equation can also be written (using \eqref{eq:rho1bis}) as
\begin{equation}\label{eq:V1}
\rho \partial_t  \mathscr V+ \rho \V\cdot\na_x \mathscr V +\rho \nabla_xP=\lambda \rho(u_{ext}-\mathscr V).
\end{equation}

\medskip

\noindent{\bf Equation for $(X(t),V(t))$.}
We recall that $V(t)$ is defined by
$$
V  (t) = 
\frac{1}{m} \ds\int_{\R^N}  j (t,x) \ud x= \frac{1}{m} \ds\int_{\R^N}  \rho( t,x) \V(t,x) \ud x$$
and integrating \eqref{eq:V1} with respect to $x$ gives
$$ V'(t) +\int_{\R^N} \rho \na_x P\ud x= \lambda \left( \frac 1 m \int \rho(t,x) u_{ext}(t,x)\ud x - V(t)\right).$$
The integral $\int_{\R^N} \rho \na_x P\ud x$ describes the macroscopic effect of the interaction potential on the  motion of the center of mass. 
Since $\rho \na_x P$ is obtained by passing to the limit in the term $\eps^{-1}\rho_\eps \na \Phi_\eps$, and 
because of the symmetry $W(-x)=W(x)$, we have
\begin{equation}\label{zeropot}
\int_{\R^N}  \eps^{-1}\rho_\eps \na \Phi_\eps \ud x = \iint_{\R^N\times \R^N} \na W(x-y) \rho_\eps (t,x) \rho_\eps(t,y)\ud x\ud y =0
\end{equation}
and so, passing to the limit $\eps\to 0$, we find
\begin{equation}\label{eq:condP}
\int_{\R^N} \rho \na_x P\ud x =0.
\end{equation}
This condition says that the interaction potential does not affect the  motion of the center of mass $X(t)$ (see also comments about this in Section~\ref{sec:rk} below)  and it yields:
$$ V'(t) = \lambda \left( \frac 1 m \int_{\R^N} \rho(t,x) u_{ext}(t,x)\ud x - V(t)\right).$$
Recalling that $\rho(t,x) = \rho_0(x-X(t))$ and
defining 
$$ g(X):=\ds\frac{1}{m} \int_{\R^N}  \rho_0(x-X) u_{ext} (t,x)\ud x
=
\frac{1}{m}  \int_{\R^N}  \rho_0(x) u_{ext} (t,x+X)\ud x,$$
we deduce that $t\mapsto (X(t),V(t))$ solves the ODE system:
\begin{equation}\label{eq:XVlim}
\begin{cases}
X' (t) = V(t), \\
V'(t) = \ds\lambda\left[g(X(t)) - V(t)\right].
\end{cases}
\end{equation}
This is a reformulation of \eqref{eq:XV1}-\eqref{eq:XV2}. Together with the initial conditions
\begin{equation}\label{eq:iniXV}
X(0)=\frac{1}{m}  \iint_{\R^N\times\R^N}  xf_0^{\mathrm{init}} (x,v)\ud v \ud x
,\qquad
V(0)=\frac{1}{m}  \iint_{\R^N\times\R^N}  vf_0^{\mathrm{init}} (x,v)\ud v \ud x
\end{equation}
this system identifies the evolution of the center of mass.
Putting everything together, we deduce that $(t,x)\mapsto (\rho(t,x), \V(t,x))$ solves \eqref{eq:asymptotic}
with $t\mapsto (X(t),V(t))$ solution of \eqref{eq:XV1}-\eqref{eq:XV2}.
\\

The existence, uniqueness and regularity of $t\mapsto (X(t),V(t))$, solution of \eqref{eq:XVlim} with initial conditions  \eqref{eq:initas} follows from the fact that the nonlinear function $g:\R^N\to \R^N$ satisfies the Lipschitz estimate:
\begin{equation}\label{eq:glip}
|g(X)-g(Y)| \leq \frac 1 m \int_{\R^N}  \rho_0(x) \left| u_{ext} (t,x+X)- u_{ext} (t,x+Y)\right|\ud x \leq \|\nabla_xu_{ext}\|_{L^\infty}|X-Y|,
\end{equation}
together with the bound $|g(X)|\leq m\|u_{ext}\|_{L^\infty}$, owing to \eqref{hyp_uext}.
In particular, the function $t\mapsto X(t)$ is defined globally in time and at least of class $C^1$. In fact we have:
\begin{lemma}\label{Lem:edo}
The system \eqref{eq:XVlim}-\eqref{eq:iniXV} admits a unique solution, defined and of class $C^1$ on $[0,\infty)$. 
\end{lemma}
This lemma makes the definition of the density
$\rho(t,x)=\rho_0(x-X(t))$ meaningful.

\subsection{Remarks and Comments} \label{sec:rk}
We note that we could forgo the characteristic system \eqref{eq:XV1}-\eqref{eq:XV2} and rewrite the asymptotic problem  
as the following system of PDEs  set on  a moving domain:
\begin{equation}\label{eq:assss}
\begin{cases}
\ds \rho(t,x) = \rho_0(x-X(t)) & \mbox{ in } \Omega(t), \\[2pt]
\ds \pa_t \rho + \na_x \cdot (\rho \V)=0 & \mbox{ in } \Omega(t), \\[2pt]
\ds \pa_t (\rho \V) + \na_x(\rho \V\otimes \V) + \rho \na_x P = \lambda \rho(u_{ext}-\V) & \mbox{ in } \Omega(t), \\[2pt]
\ds \rho \V\cdot \nu  = \rho V\cdot \nu   \mbox{ on } \pa \Omega(t) ,\qquad  \int_{\Omega(t)}\rho \na_x P\ud x=0, 
\end{cases}
\end{equation}
with $V(t)$ defined by 
$$V(t) =  \frac{1}{m} \int_{\R^N} \rho(t,x) \mathscr V(t,x)\ud x.$$

Indeed, we saw in the derivation above that the first two equations in \eqref{eq:assss} imply \eqref{eq:XV1} while integrating the momentum equation in \eqref{eq:asymptotic} (and using the boundary condition) yields
\begin{align*}
\frac{\ud}{\ud t} \int_{\Omega(t)} \rho(t,x) \V(t,x) \ud x 
& = \int_{\Omega(t)} \pa_t(\rho \V) \ud x + \int_{\pa\Omega(t)} \rho  \V V\cdot \nu \ud\sigma(x) \\
& = - \int_{\pa \Omega(t)} \rho \V \V\cdot \nu  \ud\sigma(x)- \int_{\Omega(t)}\rho \na_x P \ud x 
\\&\qquad+ \lambda \int_{\Omega(t)} \rho(u_{ext}-\V)\ud x + \int_{\pa\Omega(t)} \rho \V V\cdot \nu   \ud\sigma(x)\\
& =   - \int_{\Omega(t)}\rho \na_x P \ud x + \lambda \int_{\Omega(t)} \rho(u_{ext}-\V)\ud x,
\end{align*}
which implies   \eqref{eq:XV2} thanks to the constraint $\int_{\Omega(t)}\rho \na_x P\ud x=0$ (this constraint says  that the pressure, which enforces the divergence constraint, does not contribute to the evolution of the center of mass).

The formulation \eqref{eq:assss} of the asymptotic system is more compact, but we will prefer \eqref{eq:XV1}-\eqref{eq:XV2}-\eqref{eq:asymptotic} which shows that one can determine the evolution of $\Omega(t)$ without solving for $\V$ and is more consistent with the way solutions can be constructed.
\medskip
 
Finally, we note that the first two equations in \eqref{eq:assss} imply $-X'(t) \na \rho + \na_x\cdot (\rho \V)=0 $ in $\Omega(t)$. In particular solutions of \eqref{eq:assss} satisfy the following relation that will be useful later on:
\begin{equation}\label{eq:contweak}
\na_x \cdot (\rho(\V-V)) = 0 \quad \mbox{ in } \Omega(t), \qquad 
\rho \V\cdot \nu  = \rho V\cdot \nu   \mbox{ on } \pa \Omega(t)
\end{equation}

\begin{rmk} 
When $\lambda=0$ and $w=0$, the system  \eqref{VP1}-\eqref{VP2} can be simplified. Indeed in that case one can show that the center of mass $X_\eps(t) = \frac 1 m\int\int x  f_\eps(t,x,v)\, dx\,dv $ and velocity $V_\eps(t) = \frac 1 m\int\int v f_\eps(t,x,v)\, dx\,dv$ solves
$$
X_\eps'(t) = V_\eps(t), \qquad V_\eps'(t) = 0.$$
Up to a change of variable $x\mapsto x-X_\eps(0)-V_\eps(0)t$ and $v\mapsto v-V_\eps(0)t$, we can thus  assume that $X_\eps(t)=0$ for all $t$.
The contribution of the quadratic potential $W_a = \Lambda x\cdot x$ becomes simpler in this case since we find $W_a *\rho_\eps(t,x) = m W_a(x) + C(t)$.
We can thus write the interaction potential as
$$\Phi_\eps(t,x)=E_s*\rho_\eps(t,x)+m W_a(x)+C(t).$$
and the last term can be disregarded since it does not depend on the space variable and thus does not contribute to the force $\na_x\Phi_\eps$. 
We can thus rewrite \eqref{VP1}-\eqref{VP2} with an interaction potential that is the sum of the repulsive power potential $E_s*\rho_\eps(t,x)$ and a confining potential $m W_a(x)$:
we recover the problem that was investigated  in \cite{chiron} when $s=1$ and  in \cite{serfaty} when  $0<s\leq 1$.
\end{rmk}

\subsection{Strong Solutions for the Asymptotic Model}
\label{sec:strg}

The modulated energy method is based on introducing a suitable functional that incorporates the solution of the asymptotic system.
In order to perform the estimates necessary to prove Theorem \ref{th2}, we need to work with smooth enough solutions of this system.
For this reason we adopt the following definition:
\begin{defin}\label{def:strong}
We say that $\mathscr V(t,x)$ is a strong solution of the asymptotic system \eqref{eq:XV1}-\eqref{eq:XV2}-\eqref{eq:asymptotic} with initial condition $X^{\mathrm{init}}, V^{\mathrm{init}}, \mathscr V^{\mathrm{init}}$ (satisfying \eqref{eq:initcomp}) on the interval $[0,T]$ if $t\mapsto (X(t),V(t))$ is a $C^1$-solutions of \eqref{eq:XV1}-\eqref{eq:XV2}  and $\mathscr V\in W^{1,\infty}((0,T)\times \Omega(t))$  
solves \eqref{eq:asymptotic} for some $P\in W^{1,\infty}((0,T)\times \Omega(t))$
with initial condition \eqref{eq:initas}.
\end{defin}

It can be convenient to rewrite the asymptotic system in a  fixed domain.
Having the pair $t\mapsto (X(t),V(t))$ at hand,    the system can be recast by performing a change of variable which amounts to simply following the density/velocity field along the curve $t\mapsto (X(t),V(t))$: 
we define 
$$ \overline \V (t,x)  = \V(t,X(t)+x)-V(t), \qquad \overline P(t,x) = P(t,X(t)+x).$$
The system  \eqref{eq:asymptotic} is then equivalent to  
a system set on a \emph{fixed} domain
(this is easier to see if we use the formulation \eqref{eq:rho2}  of the continuity equation and use the velocity  equation in \eqref{eq:XVlim}):
\begin{equation}\label{eq:asymp0}
\begin{cases}
 \rho_0 \partial_t \overline \V+\rho_0 \overline  \V \cdot \na_x  \overline \V + \rho_0 \nabla_x\overline P=\lambda\rho_0 (U-\overline \V)&  \mbox{ in } (0,\infty)\times\Omega_0\\[3pt]
\na_x \cdot (\rho_0 \overline \V )=0 &  \mbox{ in } (0,\infty)\times\Omega_0, \\[3pt]
  (\rho_0 \overline \V )\cdot \nu =0 &  \mbox{ on } (0,\infty)\times\pa \Omega_0, 
\end{cases}
\end{equation}
where
\begin{align} 
U(t,x) &: =  u_{ext}(t,x+X(t))-g(X) \nonumber \\
& = \frac{1}{m} \int_{\R^N}  \rho_0(y) \left[ u_{ext}(t,x+X(t)) - u_{ext} (t,y+X(t))\right]\ud y.\label{eq:U}
\end{align}
Importantly, the velocity field $U(t,x)$ only depends on $u_{ext}(t,x)$ and on the initial conditions (since $X(t)$ is  solution of the system \eqref{eq:XVlim}) and satisfies
 $\int  \rho_0(x) U(t,x)\ud x =0$.
We recognize in system \eqref{eq:asymp0} the usual lake equation with a driving force $\lambda\rho_0 (U-\overline \V)$. 
It is supplemented with the initial condition
\begin{equation}\label{eq:init1}
 \overline \V (0,x)  = \V^{\mathrm{init}}(X(0)+x)-V(0).
\end{equation}
This formulation is natural for investigating the well-posedness of the limit system.
\\

Constructing smooth solution of \eqref{eq:asymp0} is easier when the density $\rho_0$ satisfies the non-degeneracy condition
$$ 
\rho_0(x) \geq c_0>0 \mbox{ for all } x\in \Omega_0 
$$
(which is the case when the interaction kernel is given by \eqref{kernW} with $p=2-N$ as seen in  Proposition~\ref{facile}).
We can thus obtain the existence and uniqueness of smooth solutions for \eqref{eq:asymp0}  in that case as a straightforward  extension of known results for the standard anelastic equation ($\lambda=0$) \cite{Lake1,Lake2, RT, Temam2}
and \cite[Appendix~A]{chiron}. 

The situation is much more delicate when $-N<p<2-N$ in \eqref{kernW} since  the density $\rho_0$ vanishes on $\pa\Omega_0$ (see Proposition~\ref{facile}).
Existence of solutions for the  lake equation  with such a degeneracy is investigated in  \cite{BM06,Laca,ATL}, using a vorticity-stream function reformulation which restricts the analysis to the two-dimensional case. 
These results require  $\rho_0(x)$ to be of the form $ \rho_0(x) = \varphi(x)^{1-(p+N)/2}$ with $\Omega_0 = \{\varphi>0\}$, $\varphi$ smooth, 
$\varphi|_{\partial\Omega_0}=0$ and $\na \varphi|_{\partial\Omega_0}\neq 0$, which is indeed satisfied for \eqref{kernW} in view of the explicit formula for $\rho_0$ found in \cite{CarrilloShu,CarrilloShu2d,FrMa}.
Under these assumptions, one can show the  existence and uniqueness of a solution in $C^0([0,T];W^{1,r}(\Omega_0))$, for all $r\geq 1$ (which isn't quite enough regularity for our purpose).
A development of the existence and regularity theory for this equation is clearly outside the scope of this paper and we will perform our analysis under the assumption that such a strong solution exists.
\medskip

Next, we note that when $\rho_0\geq c_0$ (that is when $p=-N+2$), then \eqref{eq:contweak} implies the 
the stronger condition  
\begin{equation}\label{bc}
 \V\cdot \nu = V(t) \cdot \nu  , \qquad \mbox{ on } \partial\Omega(t).
\end{equation}
In fact, the analysis in \cite{BM06,Laca} shows that 
this stronger condition also  holds in this degenerate case $p\in(-N,-N+2)$.
This condition \eqref{bc} plays an important role in justifying the inequality \eqref{eq:conditionsurphi}, so we give a simple justification in Lemma~\ref{Vvanishes} below.

\medskip

Finally, the definition of the modulated energy requires us to work with macroscopic quantities that are defined on $[0,T]\times \R^N$.
While the density $\rho(t,\cdot)$ is naturally defined on the whole space, the velocity $\V(t,x)$ is only defined for $x\in\Omega(t)$.
This is the final assumption made in Theorem \ref{th2}: We assume that there exists an extension, still denoted  $\mathscr V$ for the sake of simplicity, 
such that $x\mapsto \V(t,x)$ is a Lipschitz function on the whole space $\R^N$, 
  compactly supported, say in a ball that contains $\Omega(t)$ for all $0\leq t\leq T$
  and that still satisfies the condition \eqref{bc}.
The construction of such an extension is discussed   in \cite[Th.~A.1]{chiron} and in \cite[Prop.~4.1]{serfaty}, using a general result of \cite[Th.~5, VI.3]{Ste}, which apply when  $\rho_0$ is bounded from below and $\Omega_0$ is smooth enough.

\medskip

To end this section, we provide a quick proof of a simple yet important fact already mentioned above: While it is clear that the equation $\nabla\cdot(\rho U)=0$ in $\mathbb R^N$, with $\mathrm{supp}(\rho U)\subset \overline \Omega$, implies in particular the boundary condition $\rho U\cdot\nu\big|_{\partial\Omega}=0$, it is less obvious that it also implies $U\cdot\nu \big|_{\partial\Omega}=0$ when $\rho$ vanishes sublinearly on $\pa\Omega$.
This is shown in \cite{BM06,Laca}  in the context of the analysis 
of the lake equation with with degenerate densities and we provide here a simple proof at the price of further assumptions on the geometry of the domain and $\rho$, whih hold in our framework:
\begin{lemma}\label{Vvanishes}
Let $\rho: \mathbb R^N\rightarrow \mathbb R$  and $U: \mathbb R^N\rightarrow \mathbb R^N$ be continuous functions such that 
\begin{itemize}
\item $\mathrm{supp}(\rho \, U)\subset \overline \Omega$ with $\Omega$ a star-shaped open set in $\mathbb R^N$,
\item $\nabla\cdot(\rho \, U)=0$ in $\mathbb R^N$,
\item $\lim_{\mathrm{dist}(x,\partial \Omega)\to 0}
\frac{\rho(x)}{\mathrm{dist}(x,\partial \Omega)^{1-s}}=\omega>0$ for some $0<s<1$.
\end{itemize}
Then, $U\cdot\nu\big|_{\partial\Omega}=0$.
\end{lemma}

\noindent {\bf Proof.}
As mentioned above, we already know that $\rho \,U\cdot\nu\big|_{\partial\Omega}=0$. 
To prove the stronger statement, we take $0<r<1$ and note that since $\Omega$ is star-shaped, we have $r\Omega\subset \Omega$. For any $\varphi \in C^\infty_c(\mathbb R^N)$,
 we can then write:
\[
\ds\int_{\partial(r\Omega)} \rho U\cdot \nu\varphi\ud \sigma
=\int_{\mathbb R^N\setminus r\Omega} \rho U\cdot \nabla \varphi\ud x
.\]
Since $\rho U$ is continuous and supported in $ \Omega$, the right hand side can be dominated by 
\[\begin{array}{lll}
\left|\int_{\Omega\setminus r\Omega} \rho U\cdot \nabla \varphi\ud x\right|
&\leq& \|\nabla\phi\|_{L^\infty} \|\rho U\|_{L^\infty}|\Omega\setminus r\Omega|\\&\leq&
\|\nabla\phi\|_{L^\infty} \|\rho U\|_{L^\infty}|\Omega| (1- r^N) 
\end{array}\]
Since $s\in(0,1)$, we deduce:
\[\begin{array}{lll}
0&=&\ds\lim_{r\to 1}
\ds\int_{r\partial\Omega} \ds\frac{\rho(x)}{(1-r)^{1-s}} U(x)\cdot \nu(x)\varphi(x)\ud \sigma
(x)
\\[.3cm]&=&
\ds\lim_{r\to 1}
\ds\int_{\partial\Omega} \ds\frac{r^{N-1}\rho(ry)}{(1-r)^{1-s}} U(ry)\cdot \nu(ry)\varphi(ry)\ud \sigma
(y)
=\omega\ds\int_{\partial\Omega}  U(y)\cdot \nu(y)\varphi(y)\ud \sigma
(y).\end{array}\]
Since this holds for any test function $\varphi$, it follows  that $U\cdot\nu$ vanishes on $\partial\Omega$.
\QED

\section{Minimizers of the Interaction Energy}\label{sec:min}
In this Section, we recall some properties of the 
minimizer of \eqref{eq:min0} (beyond the simple cases presented in Proposition~\ref{facile}) and establish key coercivity inequalities.

The characterization of $\rho_0$ and $\Omega_0$ when $p=-N+2$ is classically related to an obstacle problem with the Laplace operator. This viewpoint extends to $p\in (-N,-N+2)$,  but involves a nonlocal fractional Laplace operator. 
We can summarize  the main features of the analysis as follows:
\begin{itemize}
\item even in the simplest case where $w=0$ in 
{(\bf H2a)}, the density $\rho_0$, which is constant in $\Omega_0$ when $p=-N+2$, is space dependent when $p\in (-N,-N+2)$ and  degenerates along $\pa\Omega_0$;
\item taking into account a perturbation $w\neq 0$  introduces additional inhomogeneities and breaks the natural radial symmetry of the problem;
\item considering $\Lambda \neq \mathbb I$ leads to interesting phenomena 
in determining the shape of  the domain $\Omega_0$. This is discussed in \cite{chiron} for the Coulombian case $p=-N+2$, and in \cite{CarrilloShu3d,CarrilloShu2d} when $p\in (-N,-N+2)$.
\end{itemize}

\subsection{Properties of the Minimizer $\rho_0$}
\label{sec:min1}
The fact that the repulsive part of the kernel $-\frac{|x|^{p}}{p}$  can be interpreted as  the fundamental solution of an integro-differential operator plays a crucial role in the proof of Theorem \ref{th2}.
More precisely, we set
$$ s:= \frac{p+N}{2}.$$
and introduce
$$
E_s(x)=
\begin{cases}
\ds\frac{\sigma(N,s)}{|x|^{N-2s}}, \qquad  &  s < \ds\frac N2,\\[.3cm]
-2 \sigma(N,s) \log|x|, & s= \ds\frac N 2,
\end{cases} \qquad \qquad \sigma(N,s) = \ds\frac{\Gamma(\frac N 2 -s)}{4^s \Gamma(s) \pi^{N/2}}.
$$  
This function is the fundamental solution of the operator  $(-\Delta)^s$, for $s\in (0,1]$:
it satisfies $(-\Delta)^s E_s =\delta$ and its Fourier transform  is given by 
$$ \widehat E_s(\xi) = \frac{1}{|\xi|^{2s}},$$
see \cite[Theorem~1.2]{Aba}.
We also introduce the homogeneous Sobolev space $\dot H^{s}(\R^N)$, for $s\in(0,1)$,
endowed with the norm
\[
\|\Phi\|_{\dot H^{s}(\R^N)}^2=\ds\int_{\R^N\times\R^N}\ds\frac{|\Phi(x)-\Phi(y)|^2}{|x-y|^{N+2s}}\ud y\ud x,
\]
which can be equivalently defined by means of the Fourier transform,
see \cite[Section~1.3, sp. Def.~1.31 \& Prop.~1.37]{BCD}.
In particular, its dual is the  functional space
\[\dot{H}^{-s}(\R^N)
=\Big\{\Phi:\R^N\to \R,\ \ds\int_{\R^N} \frac{|\widehat \Phi(\xi)|^2}{|\xi|^{2s}}\ds\frac{\ud \xi}{(2\pi)^N}<\infty\Big\}\].

With these notations, we rewrite the interaction kernel $W$ in \eqref{kernW} as (up to changing the definition of the constant $\alpha$):
\begin{equation}\label{eq:W1}
W(x) = \alpha E_s(x) + \beta \Lambda x\cdot x +w(x),
\qquad s\in(0,1],\qquad  \alpha,\beta>0.
\end{equation}

We will denote the attractive part of the potential by
\begin{equation}
W_a(x):= \Lambda x\cdot x = \sum_{j=1}^N\frac{x_j^2}{2\lambda_j^2}.
\end{equation}
Note that when $\int_{\R^N} \rho(x)\, dx =m$ and $\int_{\R^N} x \rho(x)\, dx =0$, then 
\begin{equation}\label{eq:rhoWa}
W_a *\rho(x) = m W_a(x) + C, \qquad C = \int_{\R^N} \Lambda y\cdot y\rho(y)\, dy.
\end{equation}

When  $\Lambda =\mathbb I$ and $w=0$, the existence-uniqueness statement in Proposition~\ref{facile} can be completed by an explicit expression of the minimizer $\rho_0$
and domain $\Omega_0$ (see
\cite[Theorem 3.12]{CDM} for $s=1$ and \cite{Imbert,CaVa,CarrilloShu,FrMa} for $s\in (0,1)$):
 we have
\begin{equation}\label{eq:minimizer} 
\rho_0(x) =
\begin{cases}
\ds  m N \chi_{B_{R}}(x) & \mbox{ when } s=1,\\
\ds \frac{m N}{R^{1-s}} K_{s,N} (R^2 - |x|^2)_+^{1-s} & \mbox{ when } s\in(0,1),
\end{cases}
\end{equation}
where the radius $R$ is determined by the mass constraint $\int_{\mathbb R^N} \rho_0(x)\ud x =m$. 
When $\Lambda\neq \mathbb I$, there is a competition between the radial symmetry of the repulsive kernel $E_s$ and the 
elliptic symmetry of the attractive kernel $W_a$.
In order to state the results, we introduce the following definition: 
we associate to a positive definite matrix $A=\mathrm{diag}(1/a_1^2,\dots,1/a_N^2)$,  the ellipsoids
\[
\mathcal E_A=\{x\in \R^N\,;\, A x\cdot  x \leq 1\}.\]
We then have:
\begin{proposition}[Case $\Lambda \neq \mathbb I$, $w=0$]\label{prop:ellipsoid}
Assume that $W$ is given by \eqref{eq:W1} with $s\in(0,1]$  (which corresponds to $p\in (-N,-N+2]$) and $\Lambda \neq \mathbb I$. 
Then for all $m>0$, there exists a unique minimizer $\rho_0\in \mathcal M_+^1$ such that 
$\mathscr E[\rho_0]=\mathscr E_m$ and having the center of mass at $0$. 
It satisfies:
\begin{itemize}
    \item when $s=1$ ($p=-N+2$), $\rho_0$ is in $L^\infty(\R^N)$ and there exists positive coefficients $a_1,...,a_N$ 
     such that 
    $\rho_0(x)= c_0\chi_{\mathcal  E_A}(x)$ 
    with  $A=\mathrm{diag}(1/a_1^2,\dots,1/a_N^2)$;
    \item when $s\in(0,1)$ ($-N<p < -N+2$),  then $\rho_0 \in C^{1-s}(\R^N)$. When $N=2$, there exists a matrix $A=\mathrm{diag}(1/a_1^2,1/a_2^2)$  such that $\rho_0=c_0 (1-Ax\cdot x)_+^{1-s}$.
    \end{itemize}   
\end{proposition}
The relation between the matrices $\Lambda$ and $A$ is highly non-trivial. We refer to \cite{chiron} for further discussion and references in the Coulombian case $s=1$.
When $p\in (-N,-N+2)$,  results were obtained in 
\cite{CarrilloShu2d} in dimension $2$ and extended  in \cite{Abatangelo, CarrilloShu3d} to higher dimensions when $\Lambda=\mathrm{diag}(1,..,1,1/\lambda^2)$ using the fact that the problem can be reduced to a one-dimensional problem in this case.
Some of the proofs are recalled in Section~\ref{app:ellip}, providing a practical access to the shape of the domain, and we do think that the correspondence between an ellipsoid $\mathcal E_A$ and the positive definite matrix $\Lambda$
is a general fact.
\medskip


We now turn to the general case of $W$ given by \eqref{eq:W1}, which includes the perturbation $w(x)$.
Different types of condition on $w$ have been identified that yield a unique minimizer for $\E$. Our assumption ({\bf H2b}) follows a classical approach that relies on the properties of the Fourier transform: we recast  the condition {(\bf H2b)} 
as
\begin{equation}\label{small}
 \widehat w(\xi) \geq \ds -  \frac{\kappa\alpha}{ |\xi|^{2s}} \quad \forall \xi\in \R^N , \qquad \mbox{ for some }\kappa\in (0,1)
\end{equation}
(since $w$ is even,  its Fourier transform is real valued)
and it implies that $\mathcal F[\alpha E_s + w]\geq  \frac{(1-\kappa)\alpha}{ |\xi|^{2s}}\geq 0$. This 
 implies in particular
\begin{equation}\label{eq:Econv}
\mathscr E[\mu]>0 \mbox{  for all $\mu$ signed measure with $\ds \int \ud \mu(x) =0$, $ \ds\int x  \ud \mu (x) =0$  and $\mathscr E[|\mu|]<\infty$}
\end{equation}
which is a more general condition for the uniqueness of the minimizer (it is equivalent to the strict convexity of $\lambda \mapsto \mathscr E[\lambda\rho_1+(1-\lambda)\rho_2]$ 
for any measures $\rho_1$, $\rho_2$ with same mass and center of mass).

\begin{rmk}
If we decompose $w=w_p+w_m$ with $ \widehat w_p\geq 0$ and $ \widehat w_m(\xi)\leq 0$, then  condition \eqref{small} is equivalent to $ |\widehat w_m(\xi)|\leq \ds \kappa |\xi|^{-2s}$ which can  be interpreted as a smallness condition on $(-\Delta)^{s}w_m$  in the $L^1$ norm.

We mention here another type of uniqueness result in a somewhat complementary direction: when $s=1$, $\Lambda=\mathbb I$,  and 
if $w(x)=w(|x|)$ is radially symmetric and satisfies $|\Delta w (x)|\leq \eps$ with $\eps$ small enough, then it is proved in \cite[Theorem~2.4]{ShuTadmor} that there is a unique global minimizer (which is radially symmetric).
\end{rmk}

One can show that \eqref{eq:Econv} is also satisfied when $W$ is given by \eqref{eq:W1} and $w(x) = |x|^a$ with $2\leq a\leq 4$ (see \cite{Lopes}).
However it is not clear whether the rest of our analysis could be extended to that case.
In particular, we will rely heavily on the following improvement of \eqref{eq:Econv} which follows from \eqref{small}:
\begin{lemma}
\label{l:coer}
Assume that $W$ is given by \eqref{eq:W1} where $w$ is such that \eqref{eq:H1}, \eqref{compinf} and \eqref{small} hold.
Let $\mu$ be a signed  measure on $\R^N$ such that $\int \ud \mu(x)=0$, $\int x \ud \mu(x) =0$ and $\mathscr E[|\mu|]<\infty$.
Then 
$$
\mathscr E[\mu] \geq \alpha (1-\kappa)  \| \mu \|_{\dot{H} ^{-s}(\R^N)}^2.
$$
\end{lemma}

\noindent
{\bf Proof.}
For any $\mu\in L^2(\R^N)$, 
condition \eqref{small} gives
\begin{align*}
\ds
 \iint_{\R^N\times\R^N}
\big(\alpha E_s(x-y) + w(x-y)\big) \ud \mu(y) \ud \mu(x)&=\int_{\R^N} (\alpha \widehat {E_s} (\xi)+ \widehat w(\xi) ) |\widehat \mu(\xi)|^2\ds\frac{\ud \xi}{(2\pi)^N}\nonumber \\
&\geq  \ds \alpha (1-\kappa) \ds\int_{\R^N} \ds\frac{|\widehat \mu(\xi)|^2}{|\xi|^{2s}}\ds\frac{\ud \xi}{(2\pi)^N} \nonumber \\
& \geq  \ds\alpha  (1-\kappa)  \| \mu \|_{\dot{H} ^{-s}(\R^N)}^2.  
\end{align*}
By a density argument, we deduce
\begin{equation}
 \iint_{\R^N\times\R^N}
\big(\alpha E_s(x-y) + w(x-y)\big) \mu(y) \mu(x)\ud y \ud x
\geq  \ds\alpha  (1-\kappa)  \| \mu \|_{\dot{H} ^{-s}(\R^N)}^2\label{eq:coer}
\end{equation}
for any signed  measure $\mu$  on $\R^N$ (this inequality, which we will use again later on, does not require the zero mass and zero center of mass conditions).

This implies:
\begin{align*}
\mathscr E[\mu] 
& \geq 
\ds\iint_{\R^N\times\R^N}
\ds W_a (x-y)\ud\mu(x) \ud\mu(y)
+\alpha(1-\kappa)  \| \mu \|_{\dot{H} ^{-s}(\R^N)}^2
\end{align*}
It remains to write
\[
\ds W_a(x-y)=\sum_{j=1}^N \ds\frac{|x_j-y_j|^2}{2\lambda_j^2}=\ds\sum_{j=1}^N\ds\frac{ |x_j|^2}{2\lambda_j^2}+\ds\sum_{j=1}^N
\ds\frac{|y_j|^2}{2\lambda_j^2}-\ds\sum_{j=1}^N
\ds\frac{x_jy_j}{\lambda_j^2}\]
to see that when $\int_{\R^N} \ud \mu( x) =0$ and $\int_{\R^N} x \ud \mu( x)  =0$,
 we have
 $$
\iint_{\R^N\times\R^N} W_a (x-y)  \ud\mu(x) \ud\mu(y)= 0$$
and the lemma follows.
\QED

As noted above, such an inequality yields the uniqueness of the global minimizer. We thus have the following proposition:

\begin{proposition}\label{prop:min}
Under the assumptions of Lemma \ref{l:coer},
the minimization problem \eqref{eq:min0} 
has a unique solution up to translation (in particular assumption ({\bf H2b}) implies ({\bf H1b})). 
We denote by $\rho_0$ the unique minimizer with center of mass at $0$, by $\overline{\Omega_0}$ its compact support  and by $\Phi_0=W * \rho_0$ the corresponding potential. We then have 
\begin{equation}\label{eq:rhophi}
\rho_0 \na_x \Phi_0 = 0 \mbox{ a.e.  in } \R^N
\end{equation}
and there exists a constant $A_0 = \frac{2\mathscr E[\rho_0]}{m}$ such that 
\begin{equation}\label{eq:rhophiA}
\Phi_0(x)\geq A_0, \quad \forall x\in \R^N \textrm{ and }  \Phi_0(x)=A_0 \quad \textrm{$\rho_0$-a.e.}
\end{equation}
We have $\rho_0(-x) = \rho_0(x)$.
Moreover, when $W$ is radially symmetric, then $\rho_0$ is radially symmetric and $\Omega_0=B_{R_0}$.
\end{proposition}
 
\medskip

\noindent {\bf Proof.}
The existence of a global minimizer can be established by the compactness of minimizing sequences, see for example \cite[Theorem~3.1]{Simione}
using techniques introduced in \cite{PLLcomp} and using  assumptions \eqref{eq:H1} and \eqref{compinf}.
The latter implies also that any global minimizer is compactly supported 
\cite[Theorem~1.4 and Remark~2.8]{Carr2}.

Uniqueness follows from the strict convexity of $s\mapsto \mathscr E[s\rho_1+(1-s)\rho_2]$.
Indeed, Lemma \ref{l:coer} implies that for any  two measures $\rho_1,\rho_2$ with same mass and center of mass, and for all  $r\in (0,1)$, we have (since $\mathscr 
E$ is quadratic):
$$
\mathscr E[r\rho_1+(1-r)\rho_2]-r\mathscr E[\rho_1]-(1-r)\mathscr E[\rho_2]
= - r (1-r)\mathscr E[\rho_1-\rho_2] <0
$$
with equality if and only if $\rho_1=\rho_2$. 

Finally, relations \eqref{eq:rhophi} is the usual Euler-Lagrange condition for the minimization of \eqref{energy}
and \eqref{eq:rhophiA} is a classical condition - we refer to \cite{Balag,CDM} for future discussion about this. 
Note that it is not obvious here that the sets  $\mathrm{supp}(\rho_0)$ and $\{\Phi_0=A_0\}$ are the same. In general, we only get an inclusion - see \cite{Ser_bk} for details and \cite{chiron} for counterexamples in a related context.
Such an equality is required for \eqref{eq:conditionsurphi} to hold and is thus implicitly assumed to be satisfied in Theorem \ref{th2}.
\QED

The regularity properties of both local and global minimizers have been studied for example in \cite{CDM}. 
In our setting, we recall the following result (see  \cite[Theorem~3.4, Theorem~3.6 \& Theorem~3.10]{CDM}):
\begin{proposition}
Under the assumptions of Proposition \ref{prop:min}, and if we assume in addition that $w$ is regular enough ($C^2_{\mathrm{loc}}(\R^N)$ when $s=1$ and $C_{\mathrm{loc}}^{3}(\R^N)$ when $s\in(0,1)$), then the global minimizer $\rho_0$ is such that 
\begin{itemize}
\item $\rho_0 \in L^\infty(\R^N)$ and $\Phi_0 = W* \rho_0\in C^{1,1}$ when $s=1$,  
\item $\rho_0 \in L^\infty(\R^N)\cap C^{(1-s)_-}(\R^N)$ and $\Phi_0\in C^{1,s_-}$ when $s\in (0,1)$.
\end{itemize}
Finally, for all $s\in (0,1]$ we have $\rho_0\in C^1_{\mathrm{loc}}(\Omega_0)$ 
\end{proposition}
\begin{rmk}\label{rmk:reg}
When $s\in(0,1)$, and $w=0$, we have $\rho_0$ and $\Phi_0$ in $C^{1-s}$ and $C^{1,s}$ respectively (see \cite[Remark 3.11]{CDM}). Getting this optimal regularity for general $w$ requires further regularity assumptions and isn't needed here.
\end{rmk}

Such regularity properties are obtained by using the connection between the minimization problem and an obstacle problem 
(see \cite{chiron} for the case $s=1$ and \cite{CaVa, CDM}): The definition of $E_s$ gives
\begin{equation}\label{eq:phirhoW}(-\Delta)^s\Phi_0=\alpha \rho_0+(-\Delta)^s(W_a+w)* \rho_0,
\end{equation}
Denoting $F(x)= (-\Delta)^s(W_a+w)* \rho_0(x)$,  \eqref{eq:rhophiA} implies that $\Phi_0$ solves the obstacle problem with obstacle $A_0$ and source term $F(x)$:
\begin{equation}\label{eq:obstacle}
\begin{array}{c}
\min\left\{
\Phi_0(x) - A_0, (-\Delta)^s\Phi_0 (x)-  F(x)
\right\} =0 \qquad \mbox{ in } \R^N
\end{array}
\end{equation}
(since $\rho_0\geq 0$ in $\R^N$ and $\phi_0(x)>0$ implies $\rho_0(x)=0$).
This remark is not particularly useful to determine $\rho_0$ since $F$ depends on $\rho_0$. But it is useful for proving regularity results by bootstrapping (as in \cite{CDM}) since 
$F$ is more regular than $\rho_0$ thanks to the convolution.

When $w=0$, the quadratic nature of the potential simplifies the problem since $F(x) = mW_a(x) + C$ 
does not depend on $\rho_0$ (see \eqref{eq:rhoWa}).
In the case $s=1$, we get an additional information: since 
$\Delta \Phi_0=0$ a.e. in $\{\phi_0=C_0\}\supset \mathrm{supp} (\rho_0)$, \eqref{eq:phirhoW} implies $\rho_0=\Delta(mW_a+C)$, which is constant, in $\mathrm{supp} (\rho_0)$.

\subsection{Coercivity Inequalities}
The proof of Theorem \ref{th2} will require some coercivity estimates that are more precise than that of Lemma \ref{l:coer} and which we prove here:
\begin{proposition}\label{prop:coer}
Under the assumptions of Proposition \ref{prop:min}:
\item [(i)] For any $\rho(x)$ non-negative measure on $\R^N$ with mass $m>0$, center of mass $X=\frac1{m}\int_{R^N} x \ud \rho(x)$ and  finite second order moment, we have
\begin{equation}\label{eq:coercivity0}
\mathscr E[\rho] \geq m \ds\sum_{j=1}^N \ds\int_{\R^N}\ds\frac{|x_j-X_j|^2}{\lambda_j^2}\ud \rho(x)
+\alpha(1-\kappa)
 \| \rho \|_{\dot{H} ^{-s}(\R^N)}^2.
\end{equation}

\item[(ii)] Let $\overline\rho(x) = \rho_0(x-\overline X)$ be the minimizer of \eqref{energy} 
with mass $m$ and center of mass  $\overline X\in \R^N$
and  $\overline \Phi (x) =  W* \overline \rho$ be the corresponding potential.
For any  $\rho\geq 0$ non negative measure  with mass $\int_{\R^N} \ud \rho(x)=m$, center of mass
$X= \frac{1}{m} \int_{\R^N} x\ud \rho(x)$ and  finite second order moment, we have 
\begin{equation}\label{eq:coercivity}
 \mathscr E[\rho]-\mathscr E_m
 +m^2\ds\sum_{j=1}^N \ds\frac{|X_j-\overline X_j|^2}{\lambda_j^2}
 \geq  2\int_{\R^N} (\overline \Phi - A_0) \ud \rho(x) 
+\alpha (1-\kappa) \|\rho-\overline\rho\|^2_{\dot{H}^{-s}(\R^N)}
\end{equation}
where we recall that $\mathscr E_m=\mathscr E[\overline \rho] $ denotes the minimum value of the energy in \eqref{eq:min0}.
\end{proposition}

\noindent
{\bf Proof.}
To prove \eqref{eq:coercivity0}, we recall that
\[
\ds W_a(x-y)=\sum_{j=1}^N \ds\frac{|x_j-y_j|^2}{2\lambda_j^2}=\ds\sum_{j=1}^N\ds\frac{ |x_j|^2}{2\lambda_j^2}+\ds\sum_{j=1}^N
\ds\frac{|y_j|^2}{2\lambda_j^2}-\ds\sum_{j=1}^N
\ds\frac{x_jy_j}{\lambda_j^2}\]
so that  
\[\begin{array}{lll}
\ds\iint_{\R^N\times\R^N}
\ds W_a (x-y) \rho(y)\rho(x) \ud x\ud y
&=&
m \ds\sum_{j=1}^N \ds\int_{\R^N}\ds\frac{|x_j|^2}{\lambda_j^2}\rho(x) \ud x
-\ds\sum_{j=1}^N \ds\frac{X_j^2}{\lambda_j^2}
\\
&=&m \ds\sum_{j=1}^N \ds\int_{\R^N}\ds\frac{|x_j-X_j|^2}{\lambda_j^2}\rho(x) \ud x.
\end{array}\]
The conclusion thus follows from \eqref{eq:coer} (with $\mu=\rho$).
\medskip

To prove \eqref{eq:coercivity}, we first note that $\overline \Phi(x)=\Phi_0(x-\overline X)$ and
since the problem is invariant by translation, it is enough to prove the inequality when $\overline X=0$, that is with $\overline \rho =\rho_0$ and $\overline \Phi =\Phi_0$. We write:
\begin{equation}\label{expand}\begin{array}{lll}
\mathscr E[\rho]-\mathscr E_m =\mathscr E[\rho]-\mathscr E[\rho_0] &=& \ds \int_{\R^N} \Phi \ud \rho  - \Phi_0\ud \rho_0 \\[6pt]
& = &\ds \mathscr E[\rho-\rho_0] + \int_{\R^N}  \Phi_0 \, \ud \rho (x)+ \Phi \, \ud \rho_0 (x)-2 \Phi_0\, \ud \rho_0 (x).
\end{array}\end{equation}
The last term recasts as
\[\begin{array}{lll}
-2\ds\int_{\R^N}   \Phi_0\, \ud   \rho_0(x)
&=&-2\ds\int_{\R^N}  A_0\ud  \rho_0(y)
= -2A_0\ds\int_{\R^N} \ud  \rho(y)
\end{array}\]
where we have  used the fact that $\Phi(y)=A_0$ on $\mathrm{supp}(\rho_0)$ and that 
$\int_{\mathbb R^N}\ud \rho_0(y)=m=\int_{\mathbb R^N}\ud \rho(y)$.
Next, since the kernel $ W$ is even, we get
\[
\ds\int_{\R^N}  \Phi \ud \rho_0(x)  =\ds\int_{\R^N} W* \rho\ud  \rho_0(x)=\ds\iint_{\R^N}W*  \rho_0 \ud \rho(x)
=\ds\int_{\R^N}  \Phi_0\ud \rho(x).
\]
Therefore the last three terms in \eqref{expand} become
\[
2\ds\int_{\R^N}(\underbrace{ \Phi_0-A_0)}_{\geq 0}\ud \rho(x).
\]
Let us set $\mu=\rho - \rho_0$, which is such that $\int_{\mathbb R^N} \ud \mu(x)=0$ and $\int_{\mathbb R^N} x\ud \mu(x)=m X$.
We obtain
\[\begin{array}{lll}
\ds\int_{\R^N}W_a* \mu\ud \mu(x)&=&
\ds\sum_{j=1}^N \ds\frac{1}{\lambda_j^2}
\left(\ds\iint_{\R^N\times \R^N} x_j^2\ud \mu(x)\ud \mu(y)
-\ds\iint_{\R^N\times \R^N}x_jy_j\ud\mu(x)\ud \mu(y)
\right)
\\
&=&-m^2\ds\sum_{j=1}^N \ds\frac{|X_j|^2}{\lambda_j^2}.
\end{array}\]
We conclude once again by using \eqref{eq:coer}  (with $\mu=\rho-\rho_0$) to get
\begin{align*}
\mathscr E[\rho]-\mathscr E_m  
& =  \mathscr E[\rho-\rho_0] +2 \int_{\R^N} (\Phi_0-A_0) \ud \rho(x) \\
& \geq  -m^2\ds\sum_{j=1}^N \ds\frac{|X_j|^2}{\lambda_j^2} +(1-\kappa) \| \rho-\rho_0\|^2_{\dot{H}^{-s}(\R^N)} +2\ds\int_{\R^N}(  \Phi_0-A_0) \ud \rho(x) 
\end{align*}
which completes the proof.
\QED

\section{Convergence of the Density and the Center of Mass:  Proof of Theorem \ref{th1}}\label{sec:continuity}

\subsection{Weak Solutions of the Vlasov Equation}
\label{WeakSol}

Before we investigate the behavior  of the 
weak solutions of \eqref{VP1}-\eqref{alacharge}  as $\eps\to 0$, we give
 the following statement:

\begin{proposition}\label{prop:vlasov}
Given $f^{\mathrm{init}}(x,v)$ non negative initial data in $L^1\cap L^\infty(\mathbb R^N\times \mathbb R^N)$, there exists a weak solution of \eqref{VP1} satisfying \eqref{mass_cons} and \eqref{evol_H}.
\end{proposition}

The proof proceeds from standard arguments developed for the Vlasov-Poisson system, see e.~g.~\cite{Bouchut}. We only sketch the main steps of the proof.

When the potential $\Phi:(0,\infty)\times \mathbb R^N\to \mathbb R$ is given, the Vlasov equation is a linear transport equation
which we can write in the form
\[\partial_t f+\nabla_y\cdot (Ff)=0.\]
with the variable $y=(x,v)\in \mathbb R^N\times\mathbb R^N$ and the field 
$F:(t,x,v)\mapsto (v,-\eps^{-1}\nabla\Phi(t,x)-\lambda((v-u_{ext}(t,x))$. 
Note that $\nabla_y\cdot F(t,y)=-\lambda N$ and as long as $\Phi$ is smooth enough, this equation can be solved by means of characteristics
$$\frac{\ud}{\ud t} Y(t,s,y)=F(t,Y(0,t,y)),\qquad Y(s,s,y)=y.$$
We get
$$f(t,y)=f^{\mathrm{init}}(Y(0,t,y))e^{\lambda Nt}.$$
This formula shows that $f$ is non negative and it provides 
$L^1$,  $L^\infty$ and energy estimates.
If we further assume that $f^{\mathrm{init}}$ is compactly supported, this formula also yields some estimate on the support of $f(t,\cdot)$. 

Next, we consider a sequence of mollifier $(\zeta_n)_{n\in \mathbb N}$
and the regularized convolution kernel $W_n=\zeta_n* W$.
The existence of a solution for the corresponding nonlinear equation \eqref{VP1} can be proved via a fixed point argument with the application $g\mapsto \Phi=(\zeta_n* W)* \int_{\mathbb R^N} g\ud v \mapsto f$.
Indeed, the regularity of the kernel $W_n=\zeta_n* W$ implies that the fixed point can be obtained by a direct application of the Banach theorem.

This solution $f_n$ of the nonlinear equation with the self-consistent  potential $\Phi_n=(\zeta_n* W)* \int_{\mathbb R^N} f_n \ud v$ then satisfies the same $L^1$, $L^\infty$ and energy estimates, uniformly with respect to the regularization parameter $n$.

The final step then consists in passing to the limit $n\to \infty$ in the weak formulation.
The only difficulty lies with the nonlinear term
\[
\ds\int_{0}^\infty\ds\int_{\mathbb R^N} 
\left(\ds\int_{\mathbb R^N} 
f_n\nabla_v\psi\ud v\right)  \nabla_x \Phi_n \ud x\ud t,
\]
where $\psi$ is a smooth, compactly supported test function.
Weak convergence of the potential holds in some $L^p_{\mathrm{loc}}((0,\infty)\times\mathbb R^N)$ space, thanks to standard interpolation inequalities \cite[Lemma~3,4]{Bouchut}
and convolution estimates for integrals of the type $\int_{\mathbb R^N\times\mathbb R^N}\frac{\rho(x)\rho(y)}{|x-y|^\alpha}\ud y\ud x$, see \cite[Theorem~4.3]{LL},
while strong convergence in the dual Lebesgue space of the integral $\int_{\mathbb R^N} f_n\nabla_v\psi\ud v$ is provided by applying averaging lemma techniques \cite{GLPS}.\\

\subsection{Convergence of the Density and Current}
The natural mass and energy bounds satisfied by weak solutions of \eqref{VP1} (see Proposition \ref{prop:vlasov}) imply the following uniform estimates, at the basis of the study of the asymptotic regime.

\begin{proposition}\label{prop:bounds}
Assume that the initial data satisfies  \eqref{bornener}.
Then, for all $0<T<\infty$, there exists a constant $C_T$   such that
$f_\eps(t,x,v)$ satisfies
\begin{equation}\label{eq:mass}
\ds\iint_{\R^N\times \R^N}f_\eps(t,x,v)\ud v\ud x=\ds\iint_{\R^N\times \R^N}f_\eps^{\mathrm{init}}(x,v)\ud v\ud x=m.
\end{equation}
\begin{equation}\label{kin_en}
  \sup_{0 < \eps < 1} \iint_{\R^N\times\R^N} (|x|^2+ |v|^2 ) f_\eps(t,x,v) \leq C_T \qquad \forall t\in [0,T]
\end{equation}
and 
\begin{equation}\label{eq:interen}
\mathscr E[\rho_\eps(t)] - \mathscr E_m \leq C_T \eps \qquad \forall t\in [0,T], \; \forall \eps\in (0,1).
\end{equation}
\end{proposition}
\noindent{\bf Proof.}
Given $T>0$, 
the energy functional $\mathscr H_\eps(f_\eps)$ defined in \eqref{def_H}
satisfies \eqref{evol_H} where 
$$\left|\int_{\mathbb R^N} u_{ext}\rho_\eps(t,x)\ud x\right|\leq \|u_{ext}\|_{L^\infty((0,T)\times\mathbb R^N)}\|\rho_\eps(t,\cdot)\|_{L^1(\mathbb R^N)}=m\|u_{ext}\|_{L^\infty((0,T)\times\mathbb R^N)}.
$$
With \eqref{hyp_uext} and \eqref{bornener}, it implies that, for any $0<T<\infty$, there exists $C_T$ such that 
$$\sup_{0<\eps<1, 0\leq t\leq T} \mathscr H_\eps(f_\eps)(t)\leq C_T.$$
Next, we compute
\[\ds\frac{\ud}{\ud t}\ds\iint _{\R^N\times\R^N} |x|^2f_\eps\ud v\ud x
=\ds\iint _{\R^N\times\R^N} 2x\cdot v f_\eps\ud v\ud x
\leq \ds\iint _{\R^N\times\R^N} |x|^2 f_\eps\ud v\ud x+ \ds\iint _{\R^N\times\R^N} |v|^2 f_\eps\ud v\ud x.\]
We already know that the second moment in velocity is bounded on $[0,T]$, uniformly with respect to $\eps$. Hence, we conclude by using the Gr\"onwall lemma.
\QED

The bounds of Proposition \ref{prop:bounds} allow us to show that, possibly at the price of extracting subsequences, both the density and the current have weak limits:
\begin{coro}\label{coro:bounds}
Up to a subsequence, we have
\[\begin{array}{l}
f_\eps\rightharpoonup f \text{ weakly-$\star$ in $\mathcal
  M^1([0,T]\times\mathbb R^N\times\mathbb R^N)$}, \\
\rho_\eps = \ds \int f_\eps \ud v \rightharpoonup \rho \text{ weakly-$\star$ in $\mathcal
  M^1([0,T]\times\mathbb R^N)$ and in 
  $C^0([0,T];\mathcal M^1(\R^N)$-weak-$\star$)},\\
j_\eps = \ds\int_{\R^N} v f_\eps \ud v \rightharpoonup j = \ds\int_{\R^N} v f \ud v\text{ weakly-$\star$ in $\mathcal
  M^1([0,T]\times\mathbb R^N)$}.
  \end{array}\]
 Furthermore we have $\int_{\R^N} \ud \rho(t,x)=m$ and  $\rho(t,x), j(t,x)$ solve the continuity equation \eqref{lim_mc}.
 Finally, the measure $j$ is absolutely continuous with respect to $\rho$. Precisely we have $j=\rho \mathbb V$ with 
 $\int_0^T\int_{\R^N} |\mathbb V|^2\ud \rho(t,x)<\infty$.
\end{coro}
{\bf Proof.}
The asserted compactness properties  are  immediate consequences of \eqref{eq:mass} and  \eqref{kin_en}
which imply the tightness of $f_\eps(t,x,v)$ with respect to the variables $x$ and $v$.
Furthermore, \eqref{eq:mass} and  \eqref{kin_en} imply
\begin{equation}\label{eq:tightrho}
\int_{\R^N} (1+|x|^2) \rho_\eps(t,x)\ud x\leq C_T
\end{equation}
which provides the tightness of $\rho_\eps(t,x)$ with respect to $x$.
Similarly, we get a uniform bound 
on the momentum since 
\begin{equation}\label{eq:tightj}
\begin{array}{lll}
\ds\int_{\R^N}(1+|x|)|j_\eps(t,x)|\ud x&\leq&  \ds\iint_{\R^N\times\R^N}(1+|x|) |v|\ \sqrt{f_\eps}\ \sqrt {f_\eps}\ud v\ud x 
\\&\leq& \left(\ds\iint_{\R^N\times\R^N} \lvert v \rvert^2\ f_\eps\ud v\ud x\right)^{1/2} 
\left(\ds\iint_{\R^N\times\R^N} (1+|x|^2)f_\eps\ud v\ud x\right)^{1/2}.\end{array}
\end{equation}
Extracting subsequences and letting $\eps$ go to $0$ in \eqref{eq:rhoeps},  we obtain the continuity equation \eqref{lim_mc} (at least in the sense of distributions in $(0,T)\times\R^N$).
Next, proceeding as in \cite[Section~3.5]{chiron}, 
using \eqref{eq:tightrho}-\eqref{eq:tightj},
Arzela-Ascoli theorem and diagonal extraction, 
we can  show that
\[\ds\lim_{\eps\to 0}\ds\int_{\R^N}\rho_\eps(t,x)\varphi(x)\ud x=\ds\int_{\R^N}\rho(t,x)\varphi(x)\ud x
\]
holds uniformly on $[0,T]$ for any $\varphi\in C^0(\R^N)$ such that $\lim_{|x|\to\infty} \frac{ \varphi(x)}{1+|x|^2} =0$. In particular we have $\int _{\R^N} \ud \rho(t,x) =m$.
\\

Finally, we prove that $j$ is absolutely continuous with respect to $\rho$.
This is a standard argument:
We introduce the functional 
\begin{equation}\label{funcK}
\mathscr K(\lambda,\mu)=
\left\{\begin{array}{ll}
\ds\frac12\ds\int  |\mathbb V|^2\ud \lambda & \textrm{ if $\mu=\lambda \mathbb V$},
\\
+\infty & \textrm{ otherwise,}\end{array}\right.\end{equation}
where $\lambda$ is a non negative bounded measure   on $[0,T]\times \mathbb R^N$ and $\nu$
is a vector-valued bounded measure $\mu$ on $[0,T]\times \mathbb
R^N$.
This functional is convex and lower semi-continuous as can be seen from 
 the dual formula (see \cite[Prop. 3.4]{VxBr})
\[\mathscr K(\lambda,\mu) =
\sup_{\Theta}\left\{\ds\int  \Theta\cdot\ud \mu - \ds\frac12\ds\int
   |\Theta|^2\ud \lambda \right\}\] 
where the supremum is taken over
continuous functions $\Theta:[0,T]\times \mathbb R^N\rightarrow
\mathbb R^N$. 

Since Proposition~\ref{prop:bounds} gives
$$\mathscr K(\rho_\eps,j_\eps) \leq \int_0^T \int_{\R^N\times\R^N} \frac{|v|^2}{2} f_\eps(t,x,v)\ud x\ud v\ud t \leq C_T,
$$
we deduce
$$
\mathscr K (\rho,j)
\leq \liminf_{\eps\to 0} \mathscr K(\rho_\eps,j_\eps) \leq C_T <+\infty,
$$
which means that  $j=\rho \mathbb V$, with $\int_{0}^T\int_{\R^N}  |\mathbb V|^2\ud \rho(t,x)<\infty$.
\QED

\subsection{Convergence of the Center of Mass; Identification of the Limit Density}
As pointed out  in the derivation of the limit, the motion of  the center of mass plays an important role and the following statement justifies  the behavior we formally derived in Section \ref{Sec:Formal}:
\begin{proposition}\label{eq:limitrho}
The center of mass and average velocity $(X_\eps(t),V_\eps (t) )$ defined by \eqref{eq:XVepsdef}
satisfy
$$\textrm{$X_\eps \longrightarrow X$  and $V_\eps \longrightarrow V$  uniformly over $[0,T]$},$$
where  $ (X(t),V(t))$ is the unique solution of \eqref{eq:XV1}-\eqref{eq:XV2}  with initial data \eqref{eq:initas}.
Furthermore, the limit of $\rho_\eps(t,x)$ satisfies $\rho(t,x)=\rho_0(x-X(t))$.
\end{proposition}
\noindent{\bf Proof.}
First, using \eqref{eq:mass} and \eqref{kin_en}, we obtain that $X_\eps(t)$ and $V_\eps(t)$ are bounded in $\R^N$ uniformly in $\eps$ and $t\in [0,T]$.
Furthermore, equation \eqref{VP1} implies
\begin{equation}\label{eq:XVeps}
\begin{cases}
X_\eps' (t) = V_\eps(t), \\
V_\eps'(t) =\ds \frac \lambda m\left(\int_{\R^N} \rho_\eps(t,x) u_{ext}(t,x)\ud x - V_\eps(t)\right).
\end{cases}
\end{equation}
The first equation is obtained by multiplying \eqref{VP1} by $x$ and integrating with respect to $x$ and $v$, and the second equation by multiplying \eqref{VP1} by $v$ and integrating. The derivation of this second equation uses the fact that $\na_x W(x)$ is odd so that 
\eqref{zeropot} holds.
It follows from \eqref{eq:XVeps} that the derivative $X'_\eps(t)$ and $V'_\eps(t)$ are also bounded in $\R^N$ uniformly in $\eps$ and $t\in [0,T]$.
By virtue of the Arzela-Ascoli theorem, we can thus assume that 
$$ X_\eps(t) \to \widetilde X(t), \qquad V_\eps(t)\to \widetilde V(t) \quad \mbox{ as } \eps\to 0$$
uniformly on $[0,T]$.
Now,  passing to the limit $\eps\to 0$ in the first equation of \eqref{eq:XVepsdef}, we find
\begin{equation}\label{eq:wXV}
    \widetilde X(t) = \frac{1}{m}\ds \int_{\R^N } x \ud \rho (t,x).
\end{equation}
Using the lower semicontinuity of the interaction energy $\mathcal E$ with respect to the weak-$\star$ $\mathcal M^1$ convergence,  
and \eqref{eq:interen}, we deduce
$$
\mathcal E[\rho(t,\cdot)] \leq \liminf_{\eps\to 0} \mathcal E[\rho_\eps(t,\cdot)] \leq \mathcal E_m.
$$
Since $\rho(t,\cdot)\geq 0 $ and $\int_{\R^N} \ud\rho(t,x) = m$, the uniqueness (up to translation) of the minimization problem \eqref{eq:min0} implies that $\rho(t,\cdot)$ is a translation of $\rho_0$. Equation \eqref{eq:wXV} then implies  
$$ \rho(t,x) =  \rho_0(x-\widetilde X(t)).$$
Finally, passing to the limit in \eqref{eq:XVeps}, we find that 
 $t\mapsto (\widetilde X(t), \widetilde V(t))$ is solution of 
$$
\begin{cases}
\widetilde X' (t) = \widetilde V(t), \\
\widetilde V'(t) =\ds \frac \lambda m\left(\int_{\R^N}  u_{ext}(t,x)\ud\rho(t,x) - \widetilde V(t)\right)=\ds \frac \lambda m\left(\int_{\R^N} \rho_0(x-\widetilde X(t)) u_{ext}(t,x)\ud x - \widetilde V(t)\right),
\end{cases}
$$
which is exactly the system \eqref{eq:XV1}-\eqref{eq:XV2}. 
The assumption of convergence of the initial conditions \eqref{eq:initass} and the well posedness for this system (see Lemma~\ref{Lem:edo}) imply that $(\widetilde X(t), \widetilde V(t)) = (X(t),V(t))$.
The uniqueness of the limit  implies that the whole sequence (and not just subsequences) converges.
\QED

Coming back to Lemma~\ref{coro:bounds}, 
we also deduce  that $j$ is absolutely continuous with respect to $\rho_0(.-X(t))$
 and supported in $[0,\infty)\times\overline{\Omega(t)}$.
Note that we can also show that 
$$ \widetilde V (t) = 
\frac{1}{m} \ds\int_{\R^N}  \ud j(t,x)=
\frac{1}{m} \ds\int_{\R^N} \mathbb V \ud \rho(t,x),
$$
by passing to the limit in the second equation in \eqref{eq:XVepsdef} in $\mathcal D'(0,T)$.
In the cases where $\rho_0$ is a function, $j$ is absolutely continuous  with respect to the Lebesgue measure.
All these elements prove Theorem~\ref{th1}.

\begin{rmk}
   Under the stronger assumption of Theorem \ref{th2}, we can use the coercivity inequality \eqref{eq:coercivity} to get stronger convergence of the density $\rho_\eps(t)$: The bound  \eqref{eq:interen} and   inequality \eqref{eq:coercivity} with $\overline X = X_\eps$, imply
\begin{equation}\label{eq:gjh}
C_T\eps \geq  \mathscr E[\rho_\eps(t)]-\mathscr E_m
 \geq  \alpha (1-\kappa) \|\rho_\eps(t)-\rho_0(\cdot - X_\eps(t))\|^2_{\dot{H}^{-s}(\R^N)}.
\end{equation}
Furthermore, we can write
\[
\mathscr F(\rho_0(\cdot - X_\eps(t))(\xi)
=\ds\int_{\mathbb R^N} \rho_0(x - X_\eps(t))e^{-ix\cdot \xi}
=e^{-iX_\eps(t)\cdot \xi}\widehat {\rho_0}(\xi)
\]
which converges pointwise to $e^{-iX(t)\cdot \xi}\widehat {\rho_0}(\xi)=
\mathscr F(\rho_0(\cdot - X(t))(\xi)$ as $\eps\to 0$, and is dominated by $|\widehat {\rho_0}(\xi)|$.
The Lebesgue dominated convergence theorem therefore implies that
$\rho_0(\cdot - X_\eps(t))$ converges to $\rho_0(\cdot - X(t))$ in $\dot H^{-s}(\mathbb R^N)$.
Together with \eqref{eq:gjh}, this implies that
$\rho_\eps(t,x) $ converges to $\rho_0(x-X(t))$ 
in $L^\infty(0,T;\dot H^{-s}(\R^N))$.
\end{rmk}

\section{The Flux Equation: Proof of Theorem \ref{th2}}
\label{Sec:flux}

\subsection{Modulated Energy Inequality}
We now introduce the modulated energy method used to prove the convergence of the flux $j_\eps$ to $\rho \V$.
We recall that  $f_\eps(t,x,v)$ is a  weak solution of \eqref{VP1} and we consider $(X(t),V(t),\V(t,x))$ solution of the asymptotic system \eqref{eq:XV1}-\eqref{eq:XV2}-\eqref{eq:asymptotic}.
As explained in Section~\ref{sec:strg}, 
the density $\rho(t,x)=\rho_0(x-X(t))$ is naturally defined for all $x\in \R^N$ and supported in $\Omega(t)$. 
The velocity field $\V(t,x)$, solution of the lake equation,  is only defined for $x\in\Omega(t)$ and we  consider an  extension of $\V$  and $P$ (still denoted $\V$ and $P$) to $\R^N$ as detailed in Section~\ref{sec:strg}.

We then define 
\begin{equation}\label{eq:H}\begin{array}{lll}
\mathscr H_\eps(f_\eps,X,V,\mathscr V)(t) &:=& \ds\frac{1}{2} \iint_{\R^N \times\R^N}  |v-\mathscr V(t,x)|^2 f_\eps(t,x,v)\ud v\ud x
\\&& +\ds 
 \frac{1}{2\eps} \left( \mathscr E[\rho_\eps]-\mathscr E_m + 2m \ds\sum_{j=1}^N \ds\frac{|X_{\eps,j}(t)-X_j(t)|^2}{\lambda_j^2}
 +|V_\eps(t)-V(t)|^2\right),
\end{array}\end{equation}
where $X_\eps(t)$, $V_\eps(t)$  are the center of mass and averaged velocity of $f_\eps$,  defined  by \eqref{eq:XVepsdef}.

Proposition \ref{prop:coer} implies 
\begin{equation}\label{eq:entco}
\begin{array}{lll}
\mathscr H_\eps(f_\eps,X,V,\mathscr V)(t) &\geq& \ds \frac{1}{2} \iint_{\R^N \times\R^N}  |v-\mathscr V(t,x)|^2 f_\eps(t,x,v)\ud v\ud x 
\\[8pt]
&&+ \ds\frac {(1-\kappa)\alpha} {2\eps}\|\rho_\eps-\rho_0(x-X(t))\|^2_{\dot{H}^{-s}(\R^N)}
\\[8pt]
&&+ \ds\frac 1 \eps \int_{\R^N} \rho_\eps ( \Phi_0(x-X(t)) - A_0) \ud x.
\end{array}\end{equation}
We have already proved 
the convergence of $X_\eps$, $V_\eps$ and $\rho_\eps$ towards their expected limits in Proposition~\ref{eq:limitrho}.
The 
 modulated energy strengthens this behavior and provides some control on  $j_\eps$ since 
\begin{align*}
\int_{\R^N} \left| j_\eps(t,x)  - \rho_\eps(t,x)\V\right| \ud x 
& \leq m^{ 1/ 2} \left( \iint_{\R^N\times\R^N}  |v-\mathscr V(t,x)|^2 f_\eps(t,x,v)\ud v\ud x \right)^{1/2} \\
& \leq  m^{ 1/ 2} \left( 2\mathscr H_\eps(f_\eps,X,V,\mathscr V)(t)\right)^{1/2}.
\end{align*}
Our aim is thus to prove that $\lim_{\eps\to 0}\mathscr H_\eps(f_\eps,X,V,\mathscr V)(t)= 0$. The key proposition is the following inequality:
\begin{proposition}\label{Pr:estGron}
Under the assumptions of Theorem~\ref{th2}, we have
$$
\ds\frac{\ud}{\ud t}\mathscr H_\eps(f_\eps,X,V,\mathscr V)  + \lambda \int_{\R^N\times\R^N}  |v-\mathscr V|^2 f_\eps\ud v\ud x \leq 
 C \mathscr H_\eps(f_\eps,X,V,\mathscr V) +  R_\eps(t)
 $$
for some constant $C$ which depends on $\|\V\|_{L^\infty(\R^N)}$ and $\|D \V\|_{L^\infty(\R^N)}$ but independent of $\eps$, and with
\begin{equation}\label{eq:Reps}
 R_\eps(t) : = \int_{\R^N} ( \rho_\eps \mathscr V  -j_\eps) \cdot ( \partial_t \mathscr V +\mathscr V\cdot \nabla_x\mathscr V+\lambda \mathscr V- \lambda 
u_{ext})\ud x.
 \end{equation}
\end{proposition}
Before proving this statement, we note that it implies the desired convergence of $\mathscr H_\eps(f_\eps,X,V,\mathscr V)$:
\begin{coro}\label{cor:meb}
Assume that
\begin{equation}\label{eq:ci}
\ds\lim_{\eps\to 0}\mathscr H(f_\eps^{\mathrm{init}},X^{\mathrm{init}},V^{\mathrm{init}},\mathscr V^{\mathrm{init}})=0,
\end{equation}
then
\begin{equation}\label{cvH}
\ds\lim_{\eps\to 0} \mathscr H_\eps(f_\eps,X,V,\mathscr V)(t)=0 , \qquad \mbox{ for any $0\leq t\leq T$.}
\end{equation}
\end{coro}

\noindent
{\bf Proof.}
Proposition \ref{Pr:estGron} together with Gr\"onwall lemma gives
\[
\mathscr H_\eps(f_\eps,X,V,\mathscr V) (t)
\leq e^{CT}\left(\mathscr H_\eps(f_\eps,X,V,\mathscr V)(0)+\ds\int_0^T R_\eps(s)\ud s\right) , \qquad \mbox{ for any $0\leq t\leq T$}. \]
In view of \eqref{eq:ci}, the result will follow once we show
that
\begin{equation}\label{limReps0}\ds\lim_{\eps\to 0} \int_0^T R_\eps(t) \ud t = 0.\end{equation}
We have already proved that $\rho_\eps(t,x)$ converges weakly to $\rho(t,x)=\rho_0(x-X(t))$ and $j_\eps(t,x)$ converges  weakly to $j(t,x)$, see Corollary~\ref{coro:bounds} and Proposition~\ref{eq:limitrho}. We thus have
\begin{equation}\label{limReps}
\lim_{\eps\to0}  \int_0^T R_\eps(t) \ud t
= \int_0^T\int_{\R^N} ( \rho \mathscr V  -j) \cdot ( \partial_t \mathscr V +\mathscr V\cdot \nabla_x\mathscr V+\lambda \mathscr V- \lambda 
u_{ext})\ud x\ud t.
\end{equation}
Since both $\rho(t,\cdot)$ and 
$j(t,\cdot)$ are supported in $\Omega(t)$, we can use 
\eqref{eq:asymptotic} (which holds in $[0,\infty)\times \overline{\Omega(t)}$) to recast \eqref{limReps} as
\begin{equation*} 
\lim_{\eps\to0}  \int_0^T R_\eps(t) \ud t
 =\int_0^T \int_{\R^N} ( \rho  \mathscr V  -j) \cdot \nabla_x P \ud x\ud t
 .
\end{equation*}
Using the continuity equation \eqref{lim_mc} we can write
\[\mathrm{div}_x j  = \pa_t \rho  = - \na_x\rho_0(x-X(t)) \cdot X'(t)  = - \na_x\rho \cdot V (t)  = -\mathrm{div}_x(\rho V(t)) \quad\mbox{ in } \mathcal D'((0,T)\times \R^N).
\]
We deduce that
\begin{align*}
\lim_{\eps\to0}  \int_0^T R_\eps(t) \ud t
 = \int_0^T \int_{\R^N} ( \rho  \mathscr V  -j) \cdot \nabla_x P \ud x\ud t
& =\int_0^T \int_{\R^N} ( \rho  \mathscr V  -\rho V(t)) \cdot \nabla_x P \ud x\ud t\\
& =\int_0^T \int_{\Omega(t)} ( \rho  \mathscr V  -\rho V(t)) \cdot \nabla_x P \ud x\ud t
\end{align*}
which vanishes thanks to \eqref{eq:contweak} thus giving \eqref{limReps0}.
\QED

\noindent {\bf Proof of Theorem \ref{th2}}.
We recall that Corollary \ref{coro:bounds} gives the weak converges of $j^\eps(t,x)$ (up to a subsequence) to some $j(t,x)$ absolutely continuous with respect to $\rho$.
We can identify the limit $j(t,x)$ by using the functional $\mathscr K$ defined in \eqref{funcK} again. Indeed, we have
\[\begin{array}{lll}
\mathscr K(\rho_\eps,j_\eps-\rho_\eps \mathscr V)
&=&\ds\frac12\ds\int_0^T\ds\int_{\R^N}\ds\frac{|j_\eps(t,x)-\rho_\eps(t,x)\mathscr V(t,x)|^2}{\rho_\eps(t,x)}\ud x\ud t 
\\&\leq& \ds\frac12\ds\int_0^T\ds\iint_{\R^N\times\R^N}|v-\mathscr V(t,x)|^2f_\eps(t,x,v)\ud v\ud x\ud t 
\leq \ds\int_0^T\mathscr H_\eps(f_\eps,X,V,\mathscr V)\ud t.
\end{array}\]
The lower-semicontinuity of $\mathcal K$ and Corollary \ref{cor:meb} imply that 
$$\mathscr K(\rho,j-\rho \mathscr V)=0.$$
The formula \eqref{funcK} thus implies that $j-\rho \mathscr V=0$ $\rho$ a.e. and since $j$ is absolutely continuous with respect to $\rho$, we get $j=\rho \mathscr V$.
 \QED

\subsection{Proof of Proposition~\ref{Pr:estGron}}
The functional $\mathscr H_\eps(f_\eps,X,V,\mathscr V)$ is a perturbation of the original energy $\mathscr H_\eps[f_\eps]$, defined in \eqref{def_H}. Accordingly, 
by using successively \eqref{evol_H}, \eqref{eq:rhoeps} and \eqref{eq:jeps}
we get
\begin{align*}
\frac{\ud}{\ud t}\mathscr H_\eps(f_\eps,X,V,\mathscr V) &
= \frac{\ud}{\ud t} \mathscr H_\eps
- \frac {\ud} {\ud t} \int_{\R^N}   \mathscr V(t,x)  \cdot j_\eps(t,x)\ud x
+ \frac{\ud}{\ud t}  \int_{\R^N}  \frac{|\mathscr V(t,x)|^2}{2} \rho_\eps(t,x)\ud x
\\
& \qquad 
+\frac{m^2 }{\eps} \frac{\ud}{\ud t}  \ds\sum_{j=1}^N\ds\frac{  |X_{\eps,j}(t)-Xj(t)|^2}{\lambda_j^2}
+\frac{\ud}{\ud t} \ds\frac{|V_\eps(t)-V(t)|^2}{2\eps}
\\
& = 
-\lambda \int_{\R^N\times\R^N}  v\cdot(v-u_{ext})f_\eps\ud v \ud x
-  \int_{\R^N}    \mathscr V \cdot \partial_t j_\eps \ud x
-\int_{\R^N} j_\eps \cdot \partial_t \mathscr V\ud x \\
& \qquad +   \int_{\R^N }  \frac {\mathscr V^2} 2  \partial_t \rho_\eps\ud x
 +  \int_{\R^N } \rho_\eps \mathscr V  \cdot \partial_t \mathscr V\ud x  
 \\
 &+\ds\frac {2 m^2}{\eps} \sum_{j=1}^N \ds\frac{(X_{\eps,j}(t)-X_j(t))\cdot (V_{\eps,j}(t)-V_j(t))}{\lambda_j^2}+\ds\frac{(V_\eps(t)-V(t))\cdot (\dot{V_\eps}(t)-\dot V(t))}{\eps}.
\\
& = 
- \lambda\int_{\R^N\times\R^N}  v^2 f_\eps\ud v \ud x
+ \lambda \int_{\R^N}    u_{ext} \cdot  j_\eps \ud x
+ \int_{\R^N}    \mathscr V \cdot  \mathrm{Div}_x \mathbb P_\eps  \ud x \\
& \qquad +  \int_{\R^N}    \mathscr V \cdot \frac 1 \eps \rho_\eps\nabla_x \Phi_\eps  \ud x
- \lambda \int_{\R^N}    \mathscr V \cdot (\rho_\eps u_{ext} - j_\eps)  \ud x
\\
& \qquad +\ds\int_{\R^N} ( \rho_\eps \mathscr V  -j_\eps) \cdot \partial_t \mathscr V\ud x-  \int_{\R^N } \frac {\mathscr V^2} 2  \nabla_x\cdot j_\eps\ud x  
\\&\qquad
+ \ds\frac {2 m^2}{\eps} \sum_{j=1}^N \ds\frac{(X_{\eps,j}(t)-X_j(t))\cdot (V_{\eps,j}(t)-V_j(t))}{\lambda_j^2}
 +\ds\frac{(V_\eps(t)-V(t))\cdot (\dot{V_\eps}(t)-\dot V(t))}{\eps}.
\end{align*}
We make a dissipation term appear
\begin{align*}
\frac{\ud}{\ud t}\mathscr H_\eps(f_\eps,X,V,\mathscr V) & = 
-\lambda \int_{\R^N\times\R^N}  |v-\mathscr V|^2 f_\eps\ud v \ud x 
+ \int_{\R^N}    \mathscr V \cdot 
\Big( \mathrm{Div}_x \mathbb P_\eps + \frac 1 \eps \rho_\eps\nabla_x \Phi_\eps  
\Big)
\ud x \\
& \qquad 
+\int_{\R^N} ( \rho_\eps \mathscr V  -j_\eps) \cdot (\partial_t \mathscr V +\lambda \mathscr V- \lambda u_{ext})\ud x
 -   \int_{\R^N } \frac {\mathscr V^2} 2   \nabla_x \cdot j_\eps\ud x 
 \\
 &\qquad
 + \frac {2 m^2}{\eps} \ds\sum_{j=1}^N\ds\frac{(X_{\eps,j}(t)-X_j(t))\cdot (V_{\eps,j}(t)-V_j(t))}{\lambda^2_j}
 +\ds\frac{(V_\eps(t)-V(t))\cdot (\dot{V_\eps}(t)-\dot V(t))}{\eps}.
\end{align*}
We now define
\begin{align*}
\mathbb P_{\mathscr V,\eps}  & =\int_{\R^N} (v-\mathscr V)\otimes (v-\mathscr V) f_\eps(t,x,v) \ud v \\
& =\mathbb P_\eps - \mathscr V\otimes j_\eps-j_\eps \otimes \mathscr V 
+\rho_\eps \mathscr V\otimes \mathscr V
\end{align*}
so that
\begin{align*}
\int_{\R^N}    \mathscr V \cdot  \mathrm{Div}_x \mathbb P_\eps  \ud x
& = \int_{\R^N}    \mathscr V \cdot  \mathrm{Div}_x \mathbb P_{\mathscr V,\eps}  \ud x
- \int_{\R^N}  ( \mathscr V\otimes j_\eps+j_\eps \otimes \mathscr V 
- \rho_\eps \mathscr V\otimes \mathscr V) \cdot \nabla_x \mathscr V  \ud x\\
& = \int_{\R^N}    \mathscr V \cdot  \mathrm{Div}_x \mathbb P_{\mathscr V,\eps}  \ud x
- \int_{\R^N} \left(
j_\eps \cdot \nabla_x \frac{|\mathscr V|^2}{2} + 
(j_\eps -\rho_\eps \mathscr V ) \mathscr V \cdot\nabla_x \mathscr V \right) \ud x.
\end{align*}
Plugging this into the equality above, we finally find
\begin{align*}
\frac{\ud}{\ud t}\mathscr H_\eps(f_\eps,X,V,\mathscr V) 
& = 
- \lambda\int_{\R^N\times\R^N}  |v-\mathscr V|^2 f_\eps\ud v \ud x + \int_{\R^N}    \mathscr V \cdot \left( \mathrm{Div}_x \mathbb P_{\mathscr V,\eps} + \frac 1 \eps \rho_\eps\nabla_x \Phi_\eps \right) \ud x    \\
& \qquad
+\int_{\R^N} ( \rho_\eps \mathscr V  -j_\eps) \cdot ( \partial_t \mathscr V +\mathscr V\cdot \nabla_x\mathscr V+\lambda \mathscr V- \lambda 
u_{ext})\ud x\\
& \qquad 
+ \frac {2 m^2}{\eps} \ds\sum_{j=1}^N\ds\frac{(X_{\eps,j}(t)-X_j(t))\cdot (V_{\eps,j}(t)-V_j(t))}{\lambda^2_j}\\
&\qquad+\ds\frac{(V_\eps(t)-V(t))\cdot ({V_\eps'}(t)-V'(t))}{\eps}\\
& = 
- \lambda\int_{\R^N\times\R^N}  |v-\mathscr V|^2 f_\eps\ud v \ud x 
+ R_\eps
- \int_{\R^N}   \mathrm{D}_x\mathscr V :   \mathbb P_{\mathscr V,\eps}  \ud x\\
&\qquad  +  \int_{\R^N}    \mathscr V \cdot \frac 1 \eps \rho_\eps\nabla_x \Phi_\eps  \ud x 
  + \frac {2 m^2}{\eps} \ds\sum_{j=1}^N\ds\frac{(X_{\eps,j}(t)-X_j(t))\cdot (V_{\eps,j}(t)-V_j(t))}{\lambda^2_j}
\\&\qquad
+\ds\frac{(V_\eps(t)-V(t))\cdot ({V_\eps'}(t)- V'(t))}{\eps},
\end{align*}
with $R_\eps$ defined by \eqref{eq:Reps}.
 
To complete the proof, it remains to show that, for a certain constant $C$, that does not depend on $\eps$
(but depends on norms of $\mathscr V$ and its derivatives) we have:
\begin{align}
&\left| \int_{\R^N}   \mathrm{D}_x \mathscr V :   \mathbb P_{\mathscr V,\eps}  \ud x \right|  \leq C \mathscr H_\eps(f_\eps,X,V,\mathscr V) ,\label{eq:ineq1}\\
&\left| \frac {2 m^2}{\eps} \ds\sum_{j=1}^N\ds\frac{(X_{\eps,j}(t)-X_j(t))\cdot (V_{\eps,j}(t)-V_j(t))}{\lambda^2_j} \right|  \leq C \mathscr H_\eps(f_\eps,X,V,\mathscr V) ,\label{eq:ineq2}\\
&\ds\frac{|(V_\eps(t)-V(t))\cdot ({V_\eps'}(t)-V'(t))|}{\eps}
\leq C\mathscr H_\eps(f_\eps,X,V,\mathscr V),  \label{eq:ineq2bis}\\
&\left| \frac 1 \eps \int_{\R^N}    \mathscr V \cdot \rho_\eps\nabla_x \Phi_\eps  \ud x   \right|   \leq C \mathscr H_\eps(f_\eps,X,V,\mathscr V) .\label{eq:ineq3}
\end{align}
The inequalities  \eqref{eq:ineq1} and  \eqref{eq:ineq2} follow immediately from Cauchy-Schwarz inequality and the definition \eqref{eq:H}.
To prove   \eqref{eq:ineq2bis}, we note that \eqref{eq:XVlim} and \eqref{eq:XVeps} imply
\[
\ds\frac{\ud}{\ud t}\big(V_{\eps}(t)-V(t)\big)=-\lambda \big(V_\eps(t)-V(t)\big)
+
\lambda \ds\int_{\R^N}(\rho_\eps(t,x)-\rho(t,x))
u_{ext}(t,x)\ud x
.\]
It follows that 
\begin{align*}
\frac{|(V_\eps(t)-V(t))\cdot ({V_\eps'}(t)-V'(t))|}{\eps}
& \leq -\lambda\ds\frac{|V_\eps(t)-V(t)|^2}{\eps}
+
\ds\frac{\lambda |V_\eps(t)-V(t)|}{\eps} \left|\ds\int_{\R^N}\rho_\eps(t,x)-\rho(t,x))
u_{ext}(t,x)\ud x\right|
\\
& \leq
-\lambda\ds\frac{|V_\eps(t)-V(t)|^2}{\eps}
\\
&\qquad +\ds\frac{\lambda |V_\eps(t)-V(t)|^2}{2\eps}
+\ds\frac{\lambda}{2\eps}\|u_{ext}(t,\cdot)\|^2_{\dot H^s(\R^N)}
 \|\rho_\eps(t,\cdot)-\rho(t,\cdot)\|^2_{\dot H^{-s}(\R^N)}
\\
& \leq \ds\frac{\lambda}{2\alpha(1-\kappa)}\|u_{ext}\|^2_{L^\infty(0,T;\dot H^s(\R^N))}
\mathscr H_\eps(f_\eps,X,V,\mathscr V)
\end{align*}
 which gives   \eqref{eq:ineq2bis}.
 
 \medskip
 
 It only remains to justify the inequality \eqref{eq:ineq3}. This is the most delicate part of the proof and will occupy the rest of this section. It  requires several steps.
 
First, we recall that the limiting density $\rho(t,x)=\rho_0(x-X(t))$ and the associated potential $\Phi(t,x) = W*\rho(t,x) = \Phi_0(x-X(t))$ satisfy
$\rho\nabla_x\Phi(t,x)=0$ a.e. in $\R^N$ for all $t>0$ (see \eqref{eq:rhophi}). 
We  can thus  write 
\begin{align}
\frac 1 \eps  \int_{\R^N}    \mathscr V \cdot  \rho_\eps \nabla_x \Phi_\eps  \ud x
 & = \frac 1 \eps \int_{\R^N}  \mathscr V \cdot  (\rho_\eps -\rho)\nabla_x( \Phi_\eps -\Phi) \ud x \nonumber\\
 &\qquad +  \frac 1 \eps  \int_{\R^N}  \rho \mathscr V \cdot   \nabla_x \Phi_\eps\ud x
 +  \frac 1 \eps  \int_{\R^N}  \rho_\eps \mathscr V \cdot   \nabla_x \Phi \ud x. \label{eqjhsfd}
\end{align}
Next, we note that using \eqref{eq:contweak} and the fact that $\rho$ is supported in $\Omega(t)$
we get 
$$ \int_{\R^N}  \rho (\mathscr V -V) \cdot   \nabla_x \Phi_\eps\ud x = 
 \int_{\Omega(t)}  \rho (\mathscr V -V) \cdot   \nabla_x \Phi_\eps\ud x = 0$$
and so
\begin{equation}\label{eqbj}
\ds\frac 1 \eps  \int_{\R^N}  \rho \mathscr V \cdot   \nabla_x \Phi_\eps\ud x
=
\ds \frac {1} \eps \int_{\R^N}   \rho V(t) \cdot \nabla_x \Phi_\eps\ud x
= -\ds\frac {V(t)}{\eps} \cdot    \int_{\R^N}   \rho_\eps \nabla_x  \Phi\ud x
\end{equation}
where the  last equality is obtained by using the fact that $\na_x \Phi_\eps  =\na_x W * \rho_\eps$ (and a similar definition for $\Phi$) where  $ \nabla_xW $ is odd.
Inserting \eqref{eqbj} in \eqref{eqjhsfd}, we deduce:
\begin{equation}\label{interm}
\frac 1 \eps  \int_{\R^N}    \mathscr V \cdot  \rho_\eps \nabla_x \Phi_\eps  \ud x
=
\frac 1 \eps \int_{\R^N}  \mathscr V \cdot  (\rho_\eps -\rho)\nabla_x( \Phi_\eps -\Phi) \ud x
 +  \frac 1 \eps  \int_{\R^N}  \rho_\eps (\mathscr V -V(t)) \cdot   \nabla_x \Phi \ud x.
\end{equation}

The second term  in  the right hand side of \eqref{interm} can be bounded by using \eqref{eq:conditionsurphi}
and the coercivity inequality  \eqref{eq:coercivity}
to get
\begin{align*}
 \frac 1 \eps \left| \int_{\R^N}  \rho_\eps (\mathscr V -V(t)) \cdot   \nabla_x \Phi \ud x\right|
 & \leq  \frac C \eps \int_{\R^N}  \rho_\eps  ( \Phi-A_0) \ud x \\
&  \leq  \frac C \eps \left(
 \mathscr E[\rho_\eps]-\mathscr E_m
 +m^2\ds\sum_{j=1}^N \ds\frac{|X_{\eps,j}-  X_j|^2}{\lambda_j^2}\right) \leq C \mathscr H_\eps(f_\eps,X,V,\mathscr V)(t) .
 \end{align*}
 We thus turn our attention to the first term in the right hand side of \eqref{interm}.
It will be easier to study this term if we use the  following notations:
\[\delta\rho=\rho_\eps-\rho,\qquad \delta\Phi=\Phi_\eps-\Phi= W* \delta\rho,\qquad \Phi_1=\alpha E_s * \delta\rho.\]
We then have   $(-\Delta_x)^s\Phi_1= \alpha\delta \rho$  and 
$\Phi_1=\delta\Phi-(W_a+w)* \delta \rho$,
and we split the quantity under consideration into three terms, that will be denoted
$\textcircled{\small 1}$, $\textcircled{\small 2}$, $\textcircled{\small 3}$ respectively:
\begin{align}
\frac \alpha \eps \int_{\R^N}  \mathscr V \cdot  (\rho_\eps -\rho)\nabla_x( \Phi_\eps -\Phi) \ud x 
 = \ds\frac\alpha \eps \int_{\R^N}  \mathscr V \cdot  \delta\rho\nabla_x\delta\Phi \ud x
&=
\ds\frac 1 {\eps} \int_{\R^N}  \mathscr V \cdot  ((-\Delta_x)^s\Phi_1)\nabla_x\Phi_1 \ud x
\nonumber\\
&\quad +\ds\frac 1 {\eps} \int_{\R^N}  \mathscr V \cdot  ((-\Delta_x)^s\Phi_1)\nabla_xW_a*\delta\rho \ud x
\nonumber\\
&\quad +\ds\frac 1 {\eps}\int_{\R^N}  \mathscr V \cdot  ((-\Delta_x)^s\Phi_1)\nabla_xw*\delta \rho \ud x.\label{eq:12}
\end{align}
In order to bound these three terms, we need to distinguish the cases $s=1$ and $0<s<1$.

When $s=1$ (this case is somewhat easier), we remark that
\[
\begin{array}{lll}
\textcircled{\small 1}&=&\ds-\frac 1 {\eps} \ds\sum_{i,j=1}^N\int_{\R^N}  \mathscr V_i \partial_{jj}^2\Phi_1\partial_{i}\Phi_1  \ud x
\\
&=&
\ds\frac 1 {\eps}\ds\sum_{i,j=1}^N\int_{\R^N}  \partial_{j}\mathscr V_i \partial_{j}\Phi_1\partial_{i}\Phi_1  \ud x
+\frac 1 {\eps} \ds\sum_{i,j=1}^N\int_{\R^N}  \mathscr V_i \partial_{j}\Phi_1\partial_{i}(\partial_{j}\Phi_1)  \ud x
\\&=&
\ds\frac 1 {\eps} \int_{\R^N}  \mathrm{D}_x \mathscr V:\nabla_x \Phi_1\otimes\nabla_x\Phi_1\ud x
-\ds\frac 1 {\eps} \int_{\R^N}  \nabla_x\cdot  \mathscr V|\nabla_x \Phi_1|^2\ud x.
\end{array}\]
We thus have (using \eqref{eq:entco})
\begin{align*}
|\textcircled{\small 1}|
& \leq  \frac{1}{\eps}C\| \mathrm{D}_x \mathscr V(t)\|_{L^\infty} \|\nabla_x\Phi_1\|^2_{L^2(\R^N)}\\
& \leq \frac{1}{\eps}C\| \mathrm{D}_x \mathscr V(t)\|_{L^\infty}
\|\nabla_x 
E_1*(\rho_\eps-\rho)\|^2_{L^2(\R^N)} \\
& \leq \frac{1}{\eps}C\| \mathrm{D}_x \mathscr V(t)\|_{L^\infty} \|\rho_\eps-\rho\|^2_{\dot{H}^{-1}(\R^N)} \\
& \leq   C\| \mathrm{D}_x \mathscr V(t)\|_{L^\infty} \mathscr H_\eps(f_\eps,X,V,\mathscr V)(t).
\end{align*}
 
 When $s\in(0,1)$, a similar (but less classical) computation gives:
\begin{align*}
\textcircled{\small 1}& =\frac 1 \eps \int_{\R^N}    (-\Delta_x)^s \Phi_1\mathscr V \cdot \nabla_x\Phi_1 \ud x\\
& = \frac{c_{N,s}}{\eps} \int_{\R^N} \int_{\R^N} \frac{ \Phi_1(x)-\Phi_1(y) }{|x-y|^{N+2s}} \nabla \Phi_1(x) \cdot \mathscr V(x) \, \ud y \, \ud x \\
            &= \frac{c_{N,s}}{2\eps} \int_{\R^N} \int_{\R^N} \frac{ \nabla_x [ |\Phi_1(x)-\Phi_1(y)|^2 ] }{|x-y|^{N+2s}} \mathscr V(x) \, \ud y  \, \ud x \\
            &= -\frac{c_{N,s}}{2\eps} \int_{\R^N} \int_{\R^N} \frac{ |\Phi_1(x)-\Phi_1(y)|^2 }{|x-y|^{N+2s}} \left[ \nabla\cdot \mathscr V(x) - (N+2s) \frac{x-y}{|x-y|^2} \cdot \mathscr V(x) \right] \, \ud y \, \ud x \\
            &= -\frac{c_{N,s}}{2\eps} \int_{\R^N} \int_{\R^N} \frac{ |\Phi_1(x)-\Phi_1(y)|^2 }{|x-y|^{N+2s}} \left[ \nabla \cdot \mathscr V(x) - \ds\frac{(N+2s)}{2} \frac{(x-y) \cdot ( \mathscr V(x) - \mathscr V(y) )}{|x-y|^2} \right] \, \ud y  \, \ud x,
\end{align*}
which can be controlled as before by 
\begin{align*}
 \| D\mathscr V\|_{L^\infty} \frac{C}{\eps} \int_{\R^N} \int_{\R^N} \frac{ |\Phi_1(x)-\Phi_1(y)|^2 }{|x-y|^{N+2s}} \, \ud y  \, \ud x
& = \| D\mathscr V\|_{L^\infty} \frac{C}{\eps}  \| \Phi_1\|^2_{\dot{H}^s(\R^N)}\\
& = \| D\mathscr V\|_{L^\infty} \frac{C}{\eps}  \| E_s * \delta\rho \|^2_{\dot{H}^s(\R^N)} \\
& = \| D\mathscr V\|_{L^\infty} \frac{C}{\eps}  \| \rho_\eps-\rho \|^2_{\dot{H}^{-s}(\R^N)} 
\\
& \leq   C\| \mathrm{D}_x \mathscr V(t)\|_{L^\infty} \mathscr H_\eps(f_\eps,X,V,\mathscr V)(t).
\end{align*}
We thus have a bound on the first term in \eqref{eq:12} for all $s\in (0,1]$.

We now turn our attention to the second term in \eqref{eq:12}.
 Using the fact that $\int_{\R^N}\delta \rho\ud x=0$, we write:
\[\nabla_xW_a*\delta\rho =
-\ds\sum_{j=1}^N\ds\frac{1}{\lambda_j^2}\ds\int_{\R^N}y_j\delta\rho(y)\ud y=-\ds\sum_{j=1}^N\ds\frac{X_{\eps,j}-X_j}{\lambda_j^2}.\]
We thus have (for all $s\in(0,1]$):
\begin{align*}
\textcircled{\small 2}
& =-\frac 1 \eps \ds\sum_{j=1}^N\ds\frac{X_{\eps,j}-X_j}{\lambda_j^2} \cdot \int_{\R^N}  \mathscr V (-\Delta_x)^s\Phi_1 \ud x\\
& \leq \frac C \eps \left(\sum_{j=1}^N\frac{ |X_{\eps,j}-X_j|^2 }{\lambda_j^2} \right)^{1/2} \| \V\|_{\dot{H}^s(\R^N)} \| \Phi_1\|_{\dot{H}^s(\R^N)} \\
& \leq \frac C \eps \left(\sum_{j=1}^N\frac{ |X_{\eps,j}-X_j|^2 }{\lambda_j^2} \right)^{1/2}\| \V\|_{\dot{H}^s(\R^N)} \|\delta\rho\|_{\dot{H}^{-s}(\R^N) }\\
& \leq   C\| \mathrm{D}_x \mathscr V(t)\|_{L^\infty} \mathscr H_\eps(f_\eps,X,V,\mathscr V)(t).
\end{align*}

Finally, we consider the last term in  \eqref{eq:12}, which involves the perturbation $w$. We again  distinguish the cases $s=1$ and $0<s<1$. For $s=1$, we write:
\[
\begin{array}{lll}\textcircled{\small 3} 
&=&\ds\frac 1 \eps \int_{\R^N}  (-\Delta_x\Phi_1) \mathscr V \cdot \nabla_x w* (-\Delta_x\Phi_1) \ud x
\\
&=&\ds\frac 1 \eps \ds\sum_{i,j,k=1}^N\int_{\R^N}  \mathscr V_i \partial^2_{jj}\Phi_1\partial_i w* \partial^2_{kk}\Phi_1 \ud x
\\
&=&\ds\frac 1 \eps \ds\sum_{i,j,k=1}^N\int_{\R^N}  \mathscr V_i \partial^2_{jj}\Phi_1 \partial^2_{kk} w* \partial_i\Phi_1 \ud x
\\
&=&-\ds\frac 1 \eps \ds\sum_{i,j,k=1}^N\int_{\R^N}  \partial_j\mathscr V_i \partial_{j}\Phi_1 \partial^2_{kk} w* \partial_i\Phi_1 \ud x
\\&&
-\ds\frac 1 \eps \ds\sum_{i,j,k=1}^N\int_{\R^N}  \mathscr V_i \partial_{j}\Phi_1 (\partial_j\partial^2_{kk} w)* \partial_i\Phi_1 \ud x.
\end{array}\]
This can be dominated  by 
\[\begin{array}{l}
\ds\frac 1\eps\left(\|D_x\mathscr V\|_{L^2}\|\nabla_x\Phi_1\|_{L^2}
\|\Delta_x w* \nabla_x\Phi_1\|_{L^\infty}
+
\|\mathscr V\|_{L^2}\|\nabla_x\Phi_1\|_{L^2}
\|\nabla_x(\Delta_x w)* \nabla_x\Phi_1\|_{L^\infty}
\right)
\\[.4cm]
\leq 
\ds\frac 1\eps\left(\|D_x\mathscr V\|_{L^2}\|\nabla_x\Phi_1\|^2_{L^2}
\|\Delta_x w\|_{L^2}
+
\|\mathscr V\|_{L^2}\|\nabla_x\Phi_1\|^2_{L^2}
\|\nabla_x(\Delta_x w)\|_{L^2}
\right).
\end{array}
\]
This is where the regularity assumption on $w$, {(\bf H2c)},
is used: it implies that \textcircled{\small 3}  can be dominated by $\frac C\eps \|\nabla_x\Phi_1\|^2_{L^2(\R^N)}=\frac C\eps\|\rho_\eps-\rho\|^2_{\dot{H}^{-1}(\R^N)}$, hence by the modulated energy.

When $0<s<1$, we write:
\[
\begin{array}{lll}\textcircled{\small 3} 
=\ds\frac 1 \eps \int_{\R^N}  \mathscr V \cdot  (-\Delta_x)^s\Phi_1 \nabla_x w* (-\Delta_x)^s \Phi_1 \ud x \\
=
\ds\frac 1 \eps \iiiint_{\R^{4N}}\frac{\Phi_1(x)-\Phi_1(y)}{|x-y|^{N+2s}} \times
\frac{\Phi_1(z)-\Phi_1(z')}{|z-z'|^{N+2s}} 
 \mathscr V(x) \cdot\nabla_x w(x-z)
 \ud x \ud y \ud z \ud z' \\
=
\ds\frac 1 {4\eps}  \iiiint_{\R^{4N}}
\frac{\Phi_1(x)-\Phi_1(y)}{|x-y|^{N+2s}} \times
\frac{\Phi_1(z)-\Phi_1(z')}{|z-z'|^{N+2s}}  \times \mathscr Z(x,y,z,z')
 \ud x \ud y \ud z\ud z',
\end{array}\]
where we have symmetrized the integrant and set:
\[\mathscr Z(x,y,z,z')= \mathscr V(x)\cdot \nabla_x w(x-z) -  \mathscr V(y) \cdot\nabla_x w(y-z)- \mathscr V(x) \cdot\nabla_x w(x-z') +  \mathscr V(y) \cdot\nabla_x w(y-z').\]
We deduce
\begin{equation}\label{fin3}
\textcircled{\small 3} \leq \ds\frac{C}{\eps} \|\Phi_1\|_{\dot{H}^s} 
 \left(\ds\iiiint_{\R^{4N}} \frac{
|\mathscr Z(x,y,z,z')|
^2}{|x-y|^{N+2s} |z-z'|^{N+2s}}\ud x \ud y \ud z\ud z'\right)^{1/2}.
\end{equation}
Introducing $\mathscr W = \na w$ (to simplify the notations), this last integral can be written as 
\[
\ds\iiiint_{\R^{4N}} \frac{
\big|\mathscr V(x)(\mathscr W(x-z)-\mathscr W(x-z'))-\mathscr V(y)(\mathscr W(y-z)-\mathscr W(y-z'))\big|
^2}{|x-y|^{N+2s} |z-z'|^{N+2s}}\ud x \ud y \ud z\ud z',
\]
and we claim that it can be bounded by $C(\|\mathscr W\|^2_{\dot H^s}\|\mathscr V\|_{L^2}^2+\|\mathscr W\|_{\dot H^s}^2\|\nabla\mathscr V\|_{L^2}^2+\|\nabla\mathscr W\|_{\dot H^s}^2\|\mathscr V\|_{L^2}^2)$.
To show this, we bound separately the integrals over the sets $\{|x-y|\geq 1\}$ and 
$\{|x-y|< 1\}$.
For the former, we note that:
\[\begin{array}{l}
\ds\iiiint_{|x-y|\geq 1} \frac{
\big|\mathscr V(x)(\mathscr W(x-z)-\mathscr W(x-z'))\big|^2}{|x-y|^{N+2s} |z-z'|^{N+2s}}\ud x \ud y \ud z\ud z'
\\[.4cm]
=
\ds\iint_{|x-y|\geq 1} \frac{|\mathscr V(x)|^2}{|x-y|^{N+2s}}
\left(\ds\iint_{\R^{2N}}
 \frac{\big|\mathscr W(x-z)-\mathscr W(x-z')\big|^2}{ |z-z'|^{N+2s}} \ud z\ud z'\right)\ud x \ud y
\\[.4cm]
=\ds\iint_{|x-y|\geq 1} \frac{|\mathscr V(x)|^2}{|x-y|^{N+2s}}
\left(\ds\iint_{\R^{2N}}
\frac{\big|\mathscr W(\zeta)-\mathscr W(\zeta')\big|^2}{ |\zeta-\zeta'|^{N+2s}} \ud \zeta\ud \zeta'\right)\ud x \ud y
\\[.4cm]
\leq 
\|\mathscr W\|^2_{\dot H^s} \ds\int_{\R^N}|\mathscr V(x)|^2\left(\ds\int_{|x-y|\geq 1} \frac{\ud y}{|x-y|^{N+2s}}\right)
\ud x =
\|\mathscr W\|^2_{\dot H^s} \ds\int_{\R^N}|\mathscr V(x)|^2\left(\ds\int_{|h|\geq 1} \frac{\ud h}{|h|^{N+2s}}\right)
\ud x
\\[.4cm]
\leq C\|\mathscr W\|^2_{\dot H^s}\|\mathscr V\|_{L^2}^2
\end{array}\]
(since $N+2s>N$).
We then estimate similarly the integral of $\frac{|\mathscr V(y)(\mathscr W(y-z)-\mathscr W(y-z'))|
^2}{|x-y|^{N+2s} |z-z'|^{N+2s}}$ over the same set $\{|x-y|\geq 1\}$.

For the integral over the set $\{|x-y|< 1\}$, we write
\[\begin{array}{l}
\ds\iiiint_{|x-y|< 1} \frac{
\big|\mathscr V(x)(\mathscr W(x-z)-\mathscr W(x-z'))-\mathscr V(y)(\mathscr W(y-z)-\mathscr W(y-z'))\big|^2}{|x-y|^{N+2s} |z-z'|^{N+2s}}\ud x \ud y \ud z\ud z'
\\[.4cm]
\leq\ds\iiiint_{|x-y|< 1}
\left|
\ds\int_0^1
 \Big(\nabla\mathscr V(x+\theta(y-x))(\mathscr W(x+\theta(y-x)-z)-\mathscr W(x+\theta(y-x)-z'))
\right.
\\[.4cm]
\qquad \left. +
 \mathscr V(x+\theta(y-x))(\nabla\mathscr W(x+\theta(y-x)-z)-\nabla\mathscr W(x+\theta(y-x)-z'))\Big)
\ud \theta
\right|^2
\ds \frac{\ud x \ud y \ud z\ud z'}{|x-y|^{N+2s-2} |z-z'|^{N+2s}}
\\[.4cm]
\leq 
2\ds\int_0^1\ds\iint_{|x-y|< 1}|\nabla\mathscr V(x+\theta(y-x))|^2
\left(\ds\iint_{\R^{2N}}
\ds\frac{\big|\mathscr W(\zeta)-\mathscr W(\zeta')\big|^2}{|\zeta-\zeta'|^{N+2s}}\ud \zeta\ud \zeta'
\right)
\ds \frac{\ud x \ud y\ud \theta}{|x-y|^{N+2s-2}}
\\[.4cm]
+
2\ds\int_0^1\ds\iint_{|x-y|< 1}|\mathscr V(x+\theta(y-x))|^2
\left(\ds\iint_{\R^{2N}}
\ds\frac{\big|\nabla\mathscr W(\zeta)-\nabla\mathscr W(\zeta')\big|^2}{|\zeta-\zeta'|^{N+2s}}\ud \zeta\ud \zeta'
\right)
\ds \frac{\ud x \ud y\ud \theta}{|x-y|^{N+2s-2}}
\\[.4cm]
\leq2\|\mathscr W\|_{\dot H^s}^2\ds\int_{|h|< 1}
\ds\frac{\ud h}{|h|^{N+2s-2}}
\ds\int_{\R^N} |\nabla\mathscr V(x)|^2\ud x
+
2\|\nabla\mathscr W\|_{\dot H^s}^2\ds\int_{|h|< 1}
\ds\frac{\ud h}{|h|^{N+2s-2}}
\ds\int_{\R^N} |\mathscr V(x)|^2\ud x
\\[.4cm]
\leq C(\|\mathscr W\|_{\dot H^s}^2\|\nabla\mathscr V\|_{L^2}^2+\|\nabla\mathscr W\|_{\dot H^s}^2\|\mathscr V\|_{L^2}^2)
,
\end{array}\]
since $N+2s-2<N$.
Hence, going back to \eqref{fin3} and using the fact that $ \|\Phi_1\|_{\dot{H}^s} =\alpha\|\rho_\eps-\rho\|_{\dot H^{-s}}$, 
we find
\begin{align*}
\textcircled{3}
& \leq \ds\frac C\eps \|\rho_\eps-\rho\|_{\dot H^{-s}}
(\|\nabla w\|_{\dot H^s}^2\|\nabla\mathscr V\|_{L^2}^2+\|D^2 w\|_{\dot H^s}^2\|\mathscr V\|_{L^2}^2
+\|\nabla w\|_{\dot H^s}^2\|\mathscr V\|_{L^2}^2)\\
& \leq \ds 
(\|\nabla w\|_{\dot H^s}^2\|\nabla\mathscr V\|_{L^2}^2+\|D^2 w\|_{\dot H^s}^2\|\mathscr V\|_{L^2}^2
+\|\nabla w\|_{\dot H^s}^2\|\mathscr V\|_{L^2}^2)\mathscr H_\eps(f_\eps,X,V,\mathscr V)(t) 
\end{align*}
(where we used   \eqref{eq:entco}).
This completes the proof of \eqref{eq:ineq3} and of Proposition~\ref{Pr:estGron}.

\subsection{Time Compactness of the Flux}
Theorem \ref{th2} only states that $j_\eps(t,x)$ converges to $j(t,x) $ in $\mathcal M^1_+((0,T)\times \R^N)$, but it is possible,  by procedding as in \cite{chiron}, to 
establish some time compactness property for $j_\eps$. This 
ensures in particular that the initial data for the limit equation is meaningful. We recall briefly this argument here:

We introduce the set
\[
\mathscr W=\big\{ \Theta_0:[0,T]\times\mathbb R^N\rightarrow \mathbb
R^N,\ \Theta_0 \text{ of class $C^1$},\ \mathrm{supp} (\Theta_0)\subset
[0,T]\times \overline{\Omega_0}, 
\ \nabla_x\cdot (\rho_{0}\Theta_0)=0\big\}.
\]
This is a closed subspace of the Banach space $C^1$, endowed with the sup norm for 
the function and its first order derivatives. Given $\Theta_0\in \mathscr W$, we set $$
\Theta(t,x)=\Theta_0(t,x+X(t))$$
which is supported in $\overline {\Omega(t)}$.
We multiply \eqref{eq:jeps} by $\Theta$ and we get
\begin{equation}
 \label{eqmttheta}
 \begin{array}{lll}
 \ds\frac{\ud}{\ud t}\ds\int_{\R^N} j_\eps \cdot \Theta\ud x 
& =& 
 \ds \int_{\R^N} j_\eps \cdot \partial_t \Theta\ud x 
 - \ds\int_{\R^N} \Theta \cdot (\nabla_x\cdot \mathbb P_\eps) \ud x 
\\&& +\lambda  \ds\int_{\R^N} (\rho_\eps u_{ext}-j_\eps)\cdot  \Theta \ud x
 - \ds \frac{1}{\eps} \int_{\R^N} \rho_\eps \Theta \cdot \nabla_x\Phi_{\eps} \ud x.
\end{array}\end{equation}
The first three terms in the right hand side can be readily dominated by using the Cauchy-Schwarz inequality
and owing to \eqref{kin_en}, they are all bounded, uniformly with respect to $0<\eps\ll1$ and  
$0\leq t\leq T$.
We thus look at the last term:
\[\begin{array}{lll}
 \ds \frac{1}{\eps} \int_{\R^N} \rho_\eps \Theta \cdot \nabla_x\Phi_{\eps} \ud x
 &=&
 \ds \frac{1}{\eps} \int_{\R^N} (\rho_\eps -\rho)\Theta \cdot \nabla_x(\Phi_{\eps} -\Phi)\ud x
 \\&&+
 \ds \frac{1}{\eps} \int_{\R^N} \rho\Theta \cdot \nabla_x\Phi_{\eps} \ud x
 + \ds \frac{1}{\eps} \int_{\R^N} \rho_\eps\Theta \cdot \nabla_x\Phi \ud x
 \\&&
 - \ds \frac{1}{\eps} \int_{\R^N} \rho\Theta \cdot \nabla_x\Phi \ud x.
 \end{array}\]
The last two terms vanish since $\Phi(t,\cdot)$ is constant on  $\mathrm{supp}(\rho(t,\cdot))\subset \overline{\Omega(t)}$, 
and 
$\mathrm{supp}(\Theta(t,\cdot))\subset  \overline{\Omega(t)}$.
The second term vanishes too since $\Theta_0\in \mathscr W$ and integrating by parts makes $\nabla\cdot(\rho\Theta)$ appear,
with $\nabla\cdot(\rho\Theta)(t,x)=\nabla\cdot(\rho_0\Theta_0)(t,x+X(t))=0$.
It follows that
\[
\left|\ds \frac{1}{\eps} \int_{\R^N} \rho_\eps \Theta \cdot \nabla_x\Phi_{\eps} \ud x\right|
=\frac{1}{\eps} \left| \int_{\R^N} (\rho_\eps -\rho)\Theta \cdot \nabla_x(\Phi_{\eps} -\Phi)\ud x\right|
\leq C_\Theta\|(\rho_\eps-\rho)(t,\cdot)\|^2_{H^{-s}(\R^N)}.\]
Going back to \eqref{eqmttheta}, we deduce that 
\[
\ds\frac{\ud}{\ud t}\ds\int_{\R^N}  j_\eps \cdot \Theta\ud x
\text{ is bounded in $L^\infty( (0,T) )$. }
\]
Since $\mathscr W$ is separable, we can make use of  a diagonal
argument to justify that, possibly at the price of 
extracting a subsequence, 
 $\int_{\R^N} j_\eps \cdot \Theta\ud x$
converges uniformly on $[0,T]$ for any $\Theta \in \mathscr W$.
\QED

\appendix

\section{Fractional Laplacians and Ellipsoids}
\label{sec:ell}

In this section, we recall some classical (and some not so classical) computations regarding the minimizers of the interaction energy $\mathscr E$ when $s\in(0,1)$ (non-Newtonian repulsion) and $\Lambda\neq I$ (not isotropic attraction).
These computations can be found elsewhere, but are gathered here for the reader's sake. 
One of our goal with this section is to show that the assumption \eqref{eq:conditionsurphi} is satisfied in some non-trivial cases.

\subsection{Proof of Proposition \ref{prop:ellipsoid}}
\label{app:ellip}
Proposition \ref{prop:ellipsoid} 
states that the energy minimizer of $\E$ when $s\in(0,1]$ and $\Lambda \neq I$ (and with $w=0$) is supported in an ellipsoid. 
As we will see below, this is not a simple extension of the radial symmetric case since the relation between the anisotropic matrix $\Lambda$ and the ellipsoidal support of the minimizer $\Omega_0$ is far from trivial (and largely not explicit).
The case $s=1$ was treated in details in \cite{chiron}, so we focus here on the case $s\in (0,1)$ and show that the equilibrium function has the form $c_0 (1-Ax\cdot x)_+^{1-s}$ for some diagonal positive matrix $A$.

Given a symmetric positive definite  matrix $A$, we introduce the associated inner product and norm
$$\langle x,y\rangle_A=Ax\cdot y, \qquad |x|_A= \sqrt{\langle x,x\rangle_A}. 
$$
Then, we consider  the function
\begin{equation}\label{def_uA}x\in \mathbb R^N\longmapsto  u_A(x) = (1-|x|_A^2)_+^{1-s}\end{equation}
which is thus supported in the ellipsoid
$$   \mathcal E_A = \{ x\in \mathbb R^N\, ;\, |x|_A \leq 1\}.$$
Note that 
the associated mass is 
\[\begin{array}{lll}
\ds\int_{\mathbb R^N} (1-Ax\cdot x)_+^{1-s}\ud x
&=&
\ds\int_{\mathbb R^N} (1-|y|^2)_+^{1-s} \prod_{i=1}^N \lambda _i \ud y
=
|\mathbb S^{N-1} |  \ds\int_{0}^1 (1-r^2)^{1-s}r^{N-1}\ud r\prod_{i=1}^N \lambda _i
\\[.4cm]&=&
|\mathbb S^{N-1} |  \ds\frac{B(N/2,2-s)}{2}\ \ds\prod_{i=1}^N \lambda _i,
\end{array}\]
with $B(x,y)=\frac{\Gamma(x)\Gamma(y)}{\Gamma(x+y)}$ the beta function.
We also recall \cite[Theorem~1.1]{Abatangelo} that $u_A$ solves
\begin{equation}
\begin{cases}
    (-\Delta)^{1-s} u_A  = \kappa_{s,N} & \mbox{ in } \mathcal E_A,\\
    u_A = 0 & \mbox{ in } \mathbb R^N\setminus \mathcal E_A,
\end{cases} 
\end{equation}
for some constant $\kappa_{s,N}$.

We want to show that $u_A$ is the global minimizer of the interaction energy $\E$ for some appropriate choice of the matrix $A$.
More precisely, we will show:
\begin{proposition}\label{prop:A}
Assume $N=2$ and $s\in (0,1)$. Let $W= \alpha E_s + \beta \Lambda x\cdot x$ with $\alpha,\beta>0$ and $\Lambda =\mathrm{diag}(1/\lambda_1^2,1/\lambda_2^2)$ a diagonal positive definite matrix.
Then there exists a matrix $A = \mathrm{diag}(1/a_1^2,1/a_2^2)$ such that the potential $\phi_A = W *u_A$ satisfies
\begin{equation}\label{eq:frostA} 
\phi_A(x) \geq A_0\mbox{ for all } x\in\R^2, \quad \mbox{ and }\quad  \phi_A (x)= A_0 \mbox{ for all } x\in \mathcal E_A.
\end{equation}
\end{proposition}
First, we note that we have
\[\begin{array}{lll}
\phi_A (x) = W* u_A(x)
&=& \alpha E_s* u_A
+\beta \ds\int_{\mathbb R^N} 
(\Lambda x\cdot x+\Lambda y\cdot y-2\Lambda x\cdot y) u_A (y)\ud y 
\\[.4cm]&=&
\alpha E_s* u_A
+m\beta \Lambda x\cdot x+ \beta\ds\int_{\mathbb R^N} \Lambda y\cdot y u_A (y)\ud y,
\end{array}
\]
and so \eqref{eq:frostA} is equivalent to
\begin{equation}\label{eq:frost}
\text{$\alpha E_s* u_A +m\beta\Lambda x\cdot x\geq A_0$ in $\R^N \qquad $  and 
$\qquad \alpha E_s* u_A +m\beta\Lambda x\cdot x=A_0$ in $\mathcal E_A$}.
 \end{equation}
The key step is thus to show that $E_s* u_A$ is a quadratic polynomial in $\mathcal E_A$.
This can be proved (in any dimension), but the difficulty is to get a formula for this quadratic polynomial (as a function of the matrix $A$) which is explicit enough to be inverted (so the matrix $A$ such that \eqref{eq:frost} holds can be determined).
This is the main reason why we restrict ourselves to the two dimensional case. 
We will show the following:
 \begin{lemma}
Given a positive definite matrix $A= \left(\begin{array}{cc} \frac 1 {a_1^2} & 0 \\ 0 &  \frac 1 {a_2^2}\end{array}\right)$, 
we define the quadratic polynomial
\begin{align*}
Q_A(x) &:=\ds\frac{(2-s)}{\pi}
\ds\frac{2^{3-2s}\beta_s}{C_s}R_s^{3-2s}  \left( x_1^2
\ds\int_0^{2\pi} 
\ds\frac{\cos^2(\tau)}{(a_1^2 \cos^2(\tau)+a_2^2\sin^2(\tau))^{2-s}}
\ud \tau \right.
\\&\qquad\qquad\qquad\qquad\qquad\qquad\qquad\qquad +\left.
x_2^2
\int_0^{2\pi} 
\ds\frac{\sin^2(\tau)}{(a_1^2 \cos^2(\tau)+a_2^2\sin^2(\tau))^{2-s}}
\ud \tau\right)
\end{align*}
where
$$
\beta_s=B(1-s,1-s),\quad 
R_s=\left(\ds\frac{\cos((1-s)\pi)}{(1-s)(3-2s)\pi} B\Big(\ds\frac12,\ds\frac{5-2s}{2}\Big)
\right)^{-1/(4-2s)},\quad
C_s=R_s\ds\frac{\cos((1-s)\pi)}{(1-s)(3-2s)\pi}.
$$
Then  $u_A(x) = (1-|x|_A^2)^{1-s}_+$ is such that
\begin{equation}\label{eq:EsuA}
E_s*u_A(x)\;
\begin{cases}
 \;   = V_0 - m Q_A(x) & \mbox{ for all } x\in \mathcal E_A\\
\; \geq V_0-m Q_A(x) & \mbox{ for all } x\notin \mathcal E_A,
\end{cases}
\end{equation}
where $V_0$ is a positive constant and $m=\int u_A(x)\, dx$.
 \end{lemma}
 \noindent{\bf Proof.} The proof follows  (with a slightly different presentation) the arguments developed in \cite{CarrilloShu2d}.
The idea is to reduce the computation to the one-dimensional case by using the following observation:
for $q>-1$, and $x=re_\theta\in \mathbb R^2$,  with $e_\theta=(\cos(\theta),\sin(\theta))$, we have 
\[\begin{array}{lll}
\ds\int_0^{2\pi} \left|x\cdot e_\tau 
    \right|^{q}
\ud \tau 
& =& r^{q} 
\ds\int_0^{2\pi} \left|\cos(\theta)
    \cos(\tau)+\sin(\theta)\sin(\tau)\right|^{q}
\ud \tau
=|x|^{q} \int_0^{2\pi}|\cos(\tau)|^{q}\ud \tau 
=\frac{1}{{\gamma_q}} |x|^q
\end{array}\]
with $\frac 1{\gamma_q}=\int_0^{2\pi}|\cos(\tau)|^{q}\ud \tau =2B(\frac{1+q}{2},\frac12)$.
Using this equality, we can write the convolution   $E_s*\rho$ as follows
\[
\ds E_s*\rho(x)=\int_{\mathbb R^2} \ds\frac{\rho(y)}{|x-y|^{2(1-s)}}\ud y
=\gamma_{2(s-1)}
\ds\int_0^{2\pi}\left(\ds\int_{\mathbb R^2} 
|(x-y)\cdot e_\tau|^{-2(1-s)}
\rho(y)\ud y\right)\ud \tau 
.\] 
Next, for a given $\tau\in [0,2\pi]$, we can expand $x=\tilde x_1e_\tau+\tilde x_2 e_\tau^\perp$ (where $e^\perp_\tau=(-\sin(\tau),\cos(\tau))$) 
so that 
\[E_s*\rho(x)
=\gamma_{2(s-1)}\ds\int_0^{2\pi}
\ds\int_{\mathbb R}
|\tilde x_1-\tilde y_1|^{-2(1-s)}
\left(\ds\int_{\mathbb R} \rho(\tilde y_1e_\tau+\tilde y_2 e_\tau^\perp) \ud \tilde y_2\right)\ud \tilde y_1\ud \tau.
\]
We will now apply this equality to the specific density $u_A(y)=(1-|y|^2_A)_+^{1-s}.$ With the new variables, we get
\[
|y|^2_A=\langle \tilde y_1 e_\tau+\tilde y_2 e_\tau^\perp, 
\tilde y_1 e_\tau+\tilde y_2 e_\tau^\perp\rangle_A=
B\tilde y\cdot \tilde y
\]
with 
\[
B=\begin{pmatrix}
\langle e_\tau, e_\tau\rangle_A & 
\langle e_\tau, e_\tau^\perp\rangle_A\\
\langle e_\tau, e_\tau^\perp\rangle_A & \langle e_\tau^\perp, e_\tau^\perp\rangle_A
    \end{pmatrix}.
\]
For a fixed $\tilde y_1$, 
we consider the roots of the second order polynomial 
\[P(\tilde y_2)=B\tilde y\cdot \tilde y-1
= \tilde y_2^2 \langle e_\tau^\perp, e_\tau^\perp\rangle_A
+ 2\tilde y_2 \tilde y_1 \langle e_\tau, e_\tau^\perp\rangle_A
+ \tilde y_1^2 \langle e_\tau, e_\tau\rangle_A-1.
\]
For $\tilde y_1$ small enough, the discriminant 
\[\Delta = 
\tilde y_1^2 \big(
\langle e_\tau, e_\tau^\perp\rangle_A^2- \langle e_\tau, e_\tau\rangle_A \langle e_\tau^\perp, e_\tau^\perp\rangle_A
\big)
+ \langle e_\tau^\perp, e_\tau^\perp\rangle_A
=
 \langle e_\tau^\perp, e_\tau^\perp\rangle_A\left(
 1-\tilde y_1^2 \ds\frac{|e_\tau|_A^2 |e_\tau^\perp|_A ^2- \langle e_\tau, e_\tau^\perp\rangle_A^2}{ \langle e_\tau^\perp, e_\tau^\perp\rangle_A}
 \right)
\]
is positive and $P(\tilde y_2)$ remains non positive 
for $\tilde y_2\in [\lambda_-,\lambda_+]$ with
$\lambda_\pm=\frac{-\tilde y_1\langle e_\tau,e_\tau^\perp\rangle_A\pm\sqrt\Delta}{\langle e_\tau^\perp,e_\tau^\perp\rangle_A}$.
This makes sense provided 
\[
|\tilde y_1|^2\leq \ds\frac{\langle e_\tau^\perp,e_\tau^\perp\rangle_A}{\langle e_\tau, e_\tau\rangle_A \langle e_\tau^\perp, e_\tau^\perp\rangle_A-
\langle e_\tau, e_\tau^\perp\rangle_A}.
\]
By direct inspection, this condition can be cast as
\[
|\tilde y_1|^2\leq a_1^2 \cos^2(\tau)+a_2^2 \sin^2(\tau)=
|e_\tau|^2_{A^{-1}}.
\]
The same limitation can be expressed as in \cite{CarrilloShu2d}
by searching for 
\[\mu_\tau =\max\{y\cdot e_\tau, \ y=\tilde y_1e_\tau + \tilde y_2e_\tau^\perp\in \mathcal E_A\}.\]
Indeed, 
the solution of this optimization problem can be found by means of a Lagrange multiplier:
$e_\tau + \frac p2 Ay=0$, $p(|y|_A^2-1)=0$, $p\geq 0$.
The constraint is necessarily saturated: at the optimum, $p>0$ and $|y|_A=1$
which leads to $p=-2e_\tau\cdot y$, $y=-\frac 2p A^{-1}e_\tau$, 
$(y\cdot e_\tau)^2=\tilde y_1^2
= A^{-1}e_\tau\cdot e_\tau=\mu_\tau ^2=a_1^2a_2^2 |e_\tau^\perp|_A^2$.
Therefore, we are led to
\begin{align}
E_s*u_A(x)&=
\ds\int_0^{2\pi}\ds\int_{-\mu_\tau}^{\mu_\tau}
|\tilde x_1 -\tilde y_1|^{-2(1-s)}
\left(\ds\int_{\lambda_-}^{\lambda_+}
(1-B\tilde y\cdot \tilde y)^{1-s}\ud \tilde y_2
\right)\ud \tilde y_1\ud \tau
\nonumber \\
&=
\ds\int_0^{2\pi}\ds\int_{-\mu_\tau}^{\mu_\tau}
|\tilde x_1 -\tilde y_1|^{-2(1-s)} 
\left(
\ds\int_{\lambda_-}^{\lambda_+}
\big(|e^\perp_\tau|_A^2(\lambda_+-\tilde y_2)(\tilde y_2-\lambda_-)\big)^{1-s}\ud \tilde y_2
\right)
\ud \tilde y_1\ud \tau \nonumber\\
&=
\ds\int_0^{2\pi}\ds\int_{-\mu_\tau}^{\mu_\tau}
|\tilde x_1 -\tilde y_1|^{-2(1-s)} |e^\perp_\tau|_A^{2-2s}(\lambda_+-\lambda_-)^{2(1-s)+1}
\underbrace{\left(
\ds\int_0^1 
\big(t(1-t)\big)^{1-s}
\ud t
\right)}_{:=\beta_s}
\ud \tilde y_1\ud \tau
\nonumber \\
&=
2^{3-2s}\beta_s
\ds\int_0^{2\pi}
\ds\frac{1}{|e_\tau^\perp|_A}\left(
\ds\int_{-\mu_\tau}^{\mu_\tau}
|\tilde x_1 -\tilde y_1|^{-2(1-s)} \left(
1-\ds\frac{\tilde y_1^2}{\mu_\tau^2}
\right)^{3/2-s}
\ud \tilde y_1\right)\ud \tau\label{eq:jhgf}
\end{align}
where we used the change of variable $\tilde y_2=\lambda_-+t(\lambda_+-\lambda_-)$.
(In fact $\beta_s=B(1-s,1-s)$.)

We now recall (see 
\cite[Appendix~A]{CarrilloShu3d}, \cite[Theorem~2]{FrMa} and \cite{CaVa,CDM,Fr1})
that the function
\[\tilde\rho: t\in \mathbb R\longmapsto C_s(1-t^2/R_s^2)_+^{3/2-s}\]
with 
\[
R_s=\left(\ds\frac{\cos((1-s)\pi)}{(1-s)(3-2s)\pi} B\Big(\ds\frac12,\ds\frac{5-2s}{2}\Big)
\right)^{-1/(4-2s)},\qquad
C_s=R_s\ds\frac{\cos((1-s)\pi)}{(1-s)(3-2s)\pi},
\]
realizes the minimum, unique up to translation,  of the functional 
\[
\tilde{\mathscr E}[\tilde \rho]
=\ds\iint_{\mathbb R\times \mathbb R}
\big(|t-t'|^{-2(1-s)} + |t-t'|^2\big)
\tilde \rho(t)\tilde \rho(t')
\ud t'\ud t
\]
over 
the set of probability measures.
In particular, it satisfies the Frostman conditions
\begin{equation}\label{eq:frost33}
\ds\int_{\mathbb R}
\big(|t-t'|^{-2(1-s)} + |t-t'|^2\big)
\tilde \rho(t')
\ud t'\left\{
\begin{array}{ll}=
V_1\quad \text{ if $t\in [-R_s,R_s]$,}
\\
>V_1\quad \text{ if $t\notin [-R_s,R_s]$,}
\end{array}
\right.
\end{equation}
with $V_1=R_s^2(\frac{1}{2(1-s)}+\frac{1}{2(3-2s)})$.

Looking back at \eqref{eq:jhgf}, we notice that the function   $t\mapsto \tilde\rho(R_st/\mu_\tau)$ appears in the integral. Using \eqref{eq:frost33} we can write:
\[\begin{array}{l}
\ds\int_{-\mu_\tau}^{\mu_\tau}
|\tilde x_1 -\tilde y_1|^{-2(1-s)} \left(
1-\ds\frac{\tilde y_1^2}{\mu_\tau^2}
\right)^{3/2-s}
\ud \tilde y_1
\\
\hspace*{2cm}\left\{
\begin{array}{l}
=\ds\frac1C\left(\ds\frac{\mu_\tau}{R_s}\right)^{2s-1}
\left(V_1-\ds\int_{- R_s}^{R_s} |\tilde x_1R_s/\mu_\tau-t'|^2\tilde \rho(t')\ud t'\right)\quad \text{if $\tilde x_1\in [-\mu_\tau,\mu_\tau]$},
\\
>
\ds\frac1C\left(\ds\frac{\mu_\tau}{R}\right)^{2s-1}
\left(V_1-\ds\int_{-R_s}^{R_s} |\tilde x_1R_s/\mu_\tau-t'|^2\tilde \rho(t')\ud t'\right)\quad \text{if $\tilde x_1\notin [-\mu_\tau,\mu_\tau]$.}
\end{array}
\right.
\end{array}
\]
Expanding $|\tilde x_1R_s/\mu_\tau-t'|^2$
and using the symmetry of $\tilde \rho$, 
 the right hand side  can be written as 
\[\ds\frac1C\left(\ds\frac{\mu_\tau}{R_s}\right)^{2s-1}
\left(V_1-\ds\int_{-R_s}^{R_s} |t'|^2\tilde \rho(t')\ud t'
\right)
-\ds\frac{\tilde x_1^2}C\left(\ds\frac{\mu_\tau}{R_s}\right)^{2s-3}
\underbrace{\ds\int_{-R_s}^{R_s} \tilde \rho(t')\ud t'}_{=1}.\]
In particular, for $x\in \mathcal E_A$,\eqref{eq:jhgf} implies
\[
E*u_A(x)=V_0-\ds\frac{2^{3-2s}\beta_s}{C_s}R^{3-2s} \ds\int_0^{2\pi} 
\underbrace{\ds\frac{1}{|e_\tau^\perp|_A}}_{=a_1a_2/\mu_\tau}\ds\frac{\tilde x_1^2}{\mu_\tau^{3-2s}}
\ud \tau
\]
for a certain constant $V_0>0$.
Using the fact that the mass of  $u_A$, $m$, is such that  $a_1a_2=\frac{(2-s)m}{\pi}$
and the equality
$\tilde x_1=x\cdot  e_\tau =x_1\cos(\tau)+x_2\sin(\tau)$, 
we finally arrive at
\[
\begin{array}{lll}
E_s*u_A(x)&=&V_0-x_1^2
\ds\frac{(2-s)m}{\pi}
\ds\frac{2^{3-2s}\beta_s}{C_s}R_s^{3-2s} \ds\int_0^{2\pi} 
\ds\frac{\cos^2(\tau)}{(a_1^2 \cos^2(\tau)+a_2^2\sin^2(\tau))^{2-s}}
\ud \tau
\\&&- 
x_2^2
\ds\frac{(2-s)m}{\pi}
\ds\frac{2^{3-2s}\beta_s}{C_s}R_s^{3-2s} \ds\int_0^{2\pi} 
\ds\frac{\sin^2(\tau)}{(a_1^2 \cos^2(\tau)+a_2^2\sin^2(\tau))^{2-s}}
\ud \tau
\end{array}\]
for $x\in \mathcal E_A$, which is the equality \eqref{eq:EsuA}.
We obtain similarly the inequality when $x\notin \mathcal E_A$.
\QED

\noindent{\bf Proof of Proposition \ref{prop:A}}
In order to prove \ref{eq:frost}, we need to identify values for $a_1$ and $a_2$ such that the quadratic polynomial in \eqref{eq:EsuA} is $-\frac{m\beta}{\alpha}\Lambda x\cdot x +C$ (for a given $\Lambda$ and any $C$).
In order to do this, we remark that the coefficients in front of $x_1^2$ and $x_2^2$ in \eqref{eq:EsuA} are the components of the gradient 
of the function
\[\zeta:(r_1,r_2)\longmapsto 
\ds\frac{(2-s)m}{\pi}
\ds\frac{2^{3-2s}\beta_s}{C_s(1-s)}R_s^{3-2s} \ds\int_0^{2\pi} 
(r_1 \cos^2(\tau)+r_2\sin^2(\tau))^{s-1}
\ud \tau
\]
evaluated at $(r_1,r_2)=(a_1^2,a_2^2)$. 

We thus want to solve 
\begin{equation}\label{eq:zetaeq}
    \nabla\zeta(a_1^2,a_2^2) = (z_1,z_2),
\end{equation}
where 
$(z_1,z_2)=-\frac{m\beta}\alpha (1/\lambda^2_1,1/\lambda_2^2)$.
We proceed as was done in \cite{chiron} in the case $s=1$ (see Remark \ref{rmk:chiron} below):
The function $\zeta$ is strictly convex, since its Hessian matrix is proportional to 
\[
M=\ds\int_0^{2\pi}
\begin{pmatrix}
\cos^4(\tau) & \cos^2(\tau)\sin^2(\tau)
\\
\cos^2(\tau)\sin^2(\tau) & \sin^4(\tau)
\end{pmatrix}\ds\frac{\ud \tau}{(r_1 \cos^2(\tau)+r_2\sin^2(\tau))^{3-s}}
\]
which is positive definite by virtue of the Cauchy-Schwarz inequality. 
We can thus use the Legendre transform of $\zeta$ to solve \eqref{eq:zetaeq}: 
More precisely, we consider the function
\[
\ell: r \longmapsto r\cdot z -\zeta(r)
\]
for a given $z\in(-\infty,0)^2$.
Since $\zeta$ is strictly convex, continuous and satisfies $ \zeta(r) \to 0 $ as $|r|\to \infty$, 
there exists a unique $r_0\in [0,+\infty)^2$ such that
$$\ell(r_0)= \sup_{r\in [0,+\infty)^2} \ell(r).$$
We then have
$$ z = \na \zeta(r_0)$$
as desired and 
we just need to check that the components of $r_0$ do not vanish,
which follows from the fact that
$$\lim_{r_2\to 0^+ }\pa_{2} \ell(r_1,r_2)=+\infty \quad \forall r_1>0 , \quad \lim_{r_1\to 0^+}\pa_1 \ell(r_1,r_2)=+\infty \quad \forall r_2>0.$$
We can thus take $(a_1^2,a_2^2)=r_0$ which identifies the matrix $A$ and completes the proof.
This construction can be extended to higher dimension, provided $\Lambda=\mathrm{diag}(1/\lambda^2_1,1,...,1)$, see \cite{CarrilloShu3d}.

\begin{rmk}\label{rmk:chiron}
When $s=1$, and with the convention that $\frac{r^0}0=\ln (r)$ for $r>0$, we find that
$\zeta(r_1,r_2)$ is proportional to 
\[\begin{array}{lll}
\ds\int_0^{2\pi} 
\ln (r_1 \cos^2(\tau)+r_2\sin^2(\tau))
\ud \tau
&=&\ds\int_0^{2\pi} 
\ln \left(\ds\frac{r_1+r_2}{2} +
\ds\frac{r_1-r_2}{2}\cos(2\tau)\right)
\ud \tau
\\&=&
\ds\int_0^{2\pi} 
\ln \left(\ds\frac{r_1+r_2}{2} +
\ds\frac{r_1-r_2}{2}\cos(\tau)\right)
\ud \tau
\end{array}\]
which  we denote as $I(\frac{r_1+r_2}{2} ,\frac{r_1-r_2}{2} )$.
For $A>B$, a simple computation gives
\[\begin{array}{lll}
\partial_A I(A,B)&=&
\ds\int_0^{2\pi} 
\ds\frac{\ud \tau}{A+
B\cos(\tau)}
=\ds\int_{-\infty}^{+\infty} \ds\frac{2\ud t}
{(1+t^2)(A+B(1-t^2)/(!+t^2)}
\\&=&
\ds\int_{-\infty}^{+\infty} \ds\frac{2\ud t}
{(A+B+(A-B)t^2}
=\ds\frac{\pi}{\sqrt{A^2-B^2}}
\end{array}\]
where we have set $t=\tan(\tau)$.
Integrating, we deduce
\[
I(A,B)=2\pi\ln\left(\ds\frac{A+\sqrt{A^2-B^2}}{2}\right)\]
which in turns proves  that $\zeta(r_1,r_2)$ is proportional to $\ln\left(\ds\frac{\sqrt{r_1}+\sqrt{r_2}}{2}\right)
$ which is the formula given in \cite[Prop.~2.4]{chiron}, see also \cite{Furman,Kellogg}.
\end{rmk}


\subsection{Property \eqref{eq:conditionsurphi}}\label{sec:ext}
Property  \eqref{eq:conditionsurphi} is crucial to the proof of Proposition \ref{Pr:estGron} and thus the modulated energy argument. 
 In this section, we establish this property in the case $w=0$ and $s\in (0,1]$.
We thus have $W=(\alpha E_s + \beta \Lambda x\cdot x)$ with $s\in (0,1]$ and $\Phi_0 = (\alpha E_s+\beta W_a)*u_A$ - the potential associated to the minimizer $u_A$ determined in the previous section. We recall that
$$\mathrm{supp }\, u_A = \mathcal E_A = \{ \Phi_0=A_0\}.$$
Our goal is thus to show:
\begin{proposition}\label{prop:extension}
Let $\V:\R^N\to\R^N$ be a Lipschitz, compactly supported function such that 
$\mathscr V \cdot \nu=0$ on $\pa \mathcal E_A$.
There exists $C$ such that
\begin{equation}\label{encorepluscool}
| \mathscr{V}(t,x) \cdot \nabla \Phi_0(x) | \le C( \Phi_0 (x)-A_0)\qquad  \mbox{for all  $ x \in
\mathbb R^N $},
\end{equation}
where $A_0$ is the constant appearing in Proposition \ref{prop:A}
\end{proposition} 

The key ingredients in the proof are the Lipschitz regularity of $\V$ (with the   boundary condition \eqref{bc}), and the precise behavior of the potential $\Phi_0 $ near the boundary of $\mathcal E_A$.

This analysis was carried out in \cite[Lemma~3.1]{chiron} in the case $s=1$ and 
generalized in \cite[Lemma~2.3 \& Appendix~B]{serfaty} for more general potential $W$, including the full range $s\in (0,1]$. The argument presented here is a particular case of \cite{serfaty}, with a simplified presentation due to the absence of singular point in our simple geometry.


We recall that $\Phi_0$ solves an obstacle problem and the optimal regularity/non-degeneracy of $\Phi_0$ near the free boundary plays a crucial role in the proof.
The fact that the geometry is rather simple will be helpful here since it implies that every free boundary point is a regular point. In fact, we can show: 

\begin{proposition}\label{prop:regular}
There exists $\delta$, $\eta$ positive (small) such that 
\begin{equation}\label{eq:phisub} 
\Phi_0(x)-A_0 \geq \eta  {\mathrm{dist}}(x,\mathcal E_A) ^{1+s}\quad  \mbox{ for $1<|x|_A\leq 1+\delta$}.
\end{equation}
\end{proposition}
\noindent{\bf Proof.}
We note that 
 $\Delta \Phi_0 = \alpha \Delta E_s* \rho_0  + \beta mN.$
 When $s=1$, we have $ \Delta E_s*  \rho_0   = - \rho_0$ and so 
\begin{equation}\label{eq:Delta1}
\Delta \Phi_0 = \beta m N >0\qquad \mbox{ in } \R^N\setminus \mathcal E_A.
\end{equation}
When $s\in(0,1)$, we can write
$$ \Delta E_s*  \rho_0   = -(-\Delta)^{1-s}(-\Delta)^s(E_s* \rho_0) = -(-\Delta)^{1-s} \rho_0 = -\lambda (-\Delta)^{1-s} (R-|x|_A^2)_+^{1-s}
$$
and we then use the following result:
\begin{lemma}\label{lem:frac}
For all $s\in (0,1)$,
there exists a constant $\eta>0$ (depending on $R$ and $s$) such that
$$ 
 -(-\Delta)^{1-s} (1-|x|_A^2)_+^{1-s} \geq  \eta\, {\mathrm{dist}}(x,\mathscr E_{A} )^{-(1-s)} 
 $$
for $x$ in a neighborhood of $\pa \mathcal E_A$ in $\R^N \setminus \mathcal E_A$.
\end{lemma}
This lemma (which is likely classical but whose proof is provided below for the reader's sake) implies, for $s\in (0,1)$,
\begin{equation}\label{eq:Deltas}
\Delta \Phi_0 \geq \frac{\lambda \eta}{2} {\mathrm{dist}}(x,\mathcal E_A) ^{-(1-s)} \mbox{ for $1<|x|_A\leq 1+\delta$}
\end{equation}
with $\delta>0$.
We  can thus proceed as in \cite[Lemma 3.2]{chiron},
by  applying Taylor's formula,
 to show that \eqref{eq:Delta1} (when $s=1$) and \eqref{eq:Deltas} (when $s\in(0,1)$) 
imply \eqref{eq:phisub}.
\QED

It remains to get an upper bound on $|\mathscr V\cdot \na\Phi_0|$. This follows from the $C^{1,s}$ regularity of $\Phi_0$ (Remark \ref{rmk:reg}).
In fact, \eqref{eq:phisub} implies that every free boundary point $x_0\in \pa\mathcal E_A$ is a regular point, so we can use \cite[Proposition B.7]{serfaty} to get:
\begin{lemma}\label{lem:regular}
    Let $x$ be such that $1<|x|_A\leq 1+\delta$ and let $x_0$ be the point on $\pa\mathcal E_0$ closest to $x$. Then
    \begin{align}
   & |\pa _{\nu_0}\Phi_0(x)| \leq C \,{\mathrm{dist}}(x,\mathcal E_A)^{s}     ,        \label{eq:derbd1}
 \\
   & |\pa _{\tau_0}\Phi_0(x)| \leq C\, {\mathrm{dist}}(x,\mathcal E_A)^{s+1}    ,         \label{eq:derbd2}
    \end{align}
    where $\nu_0=\nu(x_0)$ and $\tau_0=\nu_0^\perp$.
\end{lemma}

\noindent{\bf Proof of Proposition \ref{prop:extension}.}
Since $ \mathscr{V} $ is compactly supported and $\Phi_0(x)-A_0>0$ in $\R^N\setminus \mathcal E_A$ (with $\Phi_0$ continuous), we only need to show that there is a constant $C$ such that 
\begin{equation}\label{eq:ineqV} |\mathscr{V}(x) \cdot \nabla \Phi_0(x) | \le C( \Phi_0(x)-A_0)
\end{equation}
for $x\in \R^N\setminus \mathcal E_A$ in a small neighborhood of $\partial \mathcal E_A$ to get the result (the inequality clearly holds in $\mathcal E_A$ since both sides of the inequality vanish there).
We thus fix $x$ such that
$1<|x|_A\leq 1+\delta$
with $\delta$ such  that the results of Proposition \ref{prop:regular} and Lemma \ref{lem:regular} hold.
We denote by $x_0$ the point on $\pa\mathcal E_A$ closest to $x$ and $\nu_0=\nu(x_0)$.
We have 
\begin{align*} 
| \mathscr{V}(x) \cdot \nabla \Phi_0(x)| 
& \leq 
 |\mathscr{V}(x) \cdot\nu_0|| \pa_{\nu_0}\Phi_0(x)| + | \mathscr{V}(x) ||\pa_{\tau_0} \Phi_0(x)|\\
& \leq 
 |\mathscr{V}(x)-\mathscr{V}(x_0) \cdot\nu_0|| \pa_{\nu_0}\Phi_0(x)| + | \mathscr{V} ||\pa_{\tau_0} \Phi_0(x)|
 \end{align*}
 where we used the boundary condition \eqref{bc}.
Lemma \ref{lem:regular} and the Lipschitz regularity of $\mathscr V$ thus imply
\begin{align*} 
| \mathscr{V}(x) \cdot \nabla \Phi_0(x)| 
& \leq C |x-x_0|\,{\mathrm{dist}}(x,\mathcal E_A)^{s}             
+C \,{\mathrm{dist}}(x,\mathcal E_A)^{s+1}             \\
&\leq C \,{\mathrm{dist}}(x,\mathcal E_A)^{s+1}  .
\end{align*}
Combining this inequality with \eqref{eq:phisub}, we deduce
$$| \mathscr{V}(x) \cdot \nabla \Phi_0(x)| \leq C (\Phi_0(x)-A_0) $$
whenever $1<|x|_A\leq 1+\delta$, and the Proposition follows.
\QED


\noindent
{\bf Proof of Lemma \ref{lem:frac}.} 
We  denote $\alpha :=1-s$. For $|x|_A>1$, we write
\begin{align*} 
 -(-\Delta)^{\alpha} (1-|x|_A^2)_+^{\alpha}
& = c_{\alpha,N}\int_{\R^N} \frac{ (1-|y|_A^2)_+^{\alpha}-(1-|x|_A^2)_+^{\alpha}}{|x-y|^{N+2\alpha}}\ud y\\
&= c_{\alpha,N}\int_{\mathcal E_A} \frac{ (1-|y|_A^2)_+^{\alpha} }{|x-y|^{N+2\alpha}}\ud y.
\end{align*}
Next, we denote by $x_0$ the projection of $x$ onto $\partial \mathcal E_A$ and $d= |x-x_0| = {\mathrm{dist}}(x,\partial \mathcal E_A)$. We then write
\begin{align*} 
 -(-\Delta)^{\alpha} (1-|x|_A^2)_+^{\alpha}
&= c_{\alpha,N}\int_{Ay\cdot y\leq 1} \frac{ (1-|y|_A^2)_+^{\alpha} }{|x-x_0+x_0-y|^{N+2\alpha}}\ud y\\
&= c_{\alpha,N} d^{-2\alpha}\int_{
A(x_0-dy)\cdot (x_0-d y)\leq 1} \frac{ (1-|x_0-dy|_A^2)_+^{\alpha} }{|\frac{x-x_0}{d}+y|^{N+2\alpha}}\ud y\\
&= c_{\alpha,N} d^{-2\alpha}\int_{A(x_0-dy)\cdot (x_0-d y)\leq 1} \frac{ (1-|x_0-dy|_A^2)_+^{\alpha} }{|e+y|^{N+2\alpha}}\ud y
\end{align*}
where we used the change of variable $\frac{x_0-y}{d}\mapsto y$ and denoted  $e: = \frac{x-x_0}{d}\in \mathbb S^{N-1}$.
Since $|x_0-dy|_A^2 =1-2d A x_0\cdot y + d^2 |y|_A^2$, we get
\begin{align*}  
-(-\Delta)^{\alpha} (1-|x|_A^2)_+^{\alpha}
& =   c_{\alpha,N} d^{-2\alpha}\int_{A(x_0-dy)\cdot (x_0-d y)\leq 1} \frac{ (2d A x_0\cdot y - d^2 |y|_A^2)_+^{\alpha} }{|e+y|^{N+2\alpha}}\ud y\\
& =   c_{\alpha,N} d^{-\alpha}\int_{2Ax_0\cdot y\geq d |y|^2_A} \frac{ (2 A x_0\cdot y - d |y|_A^2)_+^{\alpha} }{|e+y|^{N+2\alpha}}\ud y.
\end{align*}
We note that $e=\frac{A x_0}{|Ax_0|}$ (since $e$ is orthogonal to $\pa \mathcal E_A$),
so passing to the limit $d=|x-x_0|\to 0$,  
the integral in the equality above converges to 
$$ \eta_0 := \int_{Ax_0\cdot y>0} \frac{ (2A x_0\cdot y )_+^{\alpha} }{|e+y|^{N+2\alpha}}\ud y
=\int_{e\cdot y>0} \frac{ (2|A x_0| e\cdot y )_+^{\alpha} }{|e+y|^{N+2\alpha}}\ud y \in(0,+\infty)
.$$
Indeed, there are no singularity since $|e+y|\geq 1$ when $e\cdot y>0$ and the integrand is bounded by $\frac{C}{|y|^{N+\alpha}}$ and thus integrable for large $|y|$.

For $d$ small enough, we thus have 
$$-(-\Delta)^{\alpha} (1-|x|_A^2)_+^{\alpha}\geq \eta d^{-\alpha}$$
with $\eta =c_{\alpha,N} \eta_0/2$.
\QED

\section*{Acknowledgments}
T. G. acknowledges the Brin Mathematics Research Center
for its invitation which has permitted 
to realize the main steps
of this work. 
T. G. is partially supported by the CNRS International Emerging Action.
A. M. is partially supported by NSF Grant  DMS-2307342.
\bibliography{VPE}
\bibliographystyle{plain}

\end{document}